\documentclass[a4paper,11pt,twoside]{amsart}

\usepackage{amsfonts}
\usepackage{amsmath}
\usepackage{amsthm}
\usepackage{amssymb}
\usepackage[all]{xy}
\usepackage{hyperref}
\usepackage{aliascnt}
\usepackage[mathscr]{euscript}

\newtheorem{thm}{Theorem}[section]

\newaliascnt{prop}{thm}
\newtheorem{prop}[prop]{Proposition}
\aliascntresetthe{prop}

\newaliascnt{lem}{thm}

\aliascntresetthe{lem}

\newaliascnt{cor}{thm}
\newtheorem{cor}[cor]{Corollary}
\aliascntresetthe{cor}

\newaliascnt{conjecture}{thm}

\aliascntresetthe{conjecture}

\newaliascnt{qn}{thm}

\aliascntresetthe{qn}

\theoremstyle{definition}

\newaliascnt{definition}{thm}
\newtheorem{definition}[definition]{Definition}
\aliascntresetthe{definition}

\newaliascnt{remark}{thm}
\newtheorem{remark}[remark]{Remark}
\aliascntresetthe{remark}

\newaliascnt{ex}{thm}
\newtheorem{ex}[ex]{Example}
\aliascntresetthe{ex}

\newaliascnt{pb}{thm}
\newtheorem{pb}[pb]{Problem}
\aliascntresetthe{pb}

\numberwithin{equation}{section}

\newcommand{\ZZ}{\mathbb{Z}}

\newcommand{\iso}{\cong}

\newcommand{\farg}{-} 

\newcommand{\st}{\mid} 
\newcommand{\comp}{\circ} 

\newcommand{\mono}{\hookrightarrow} 
\newcommand{\epi}{\twoheadrightarrow} 
\newcommand{\kk}{\Bbbk} 
\newcommand{\Hom}{\mathrm{Hom}}

\newcommand{\cat}[1]{{\mathbf{#1}}} 
\newcommand{\opp}{^{\circ}} 
\newcommand{\lto}{\longrightarrow}
\newcommand{\fun}[1]{\mathsf{#1}} 
\newcommand{\Ho}[1]{\mathrm{Ho}(#1)} 

\newcommand{\sh}[2][1]{#2[#1]} 

\newcommand{\colimm}[1]{\underset{#1}{\mathrm{colim}}\,}
\newcommand{\ho}[1]{\mathrm{Ho}(#1)} 
\newcommand{\dgC}{\scat{C}_{\scat{dg}}} 
\newcommand{\dgCacy}{\scat{C}_{\scat{dg},\text{acy}}} 
\newcommand{\dgYon}[1][\cc]{\fun{Y}^{#1}_{\mathrm{dg}}} 
\newcommand{\hproj}[1]{\mathrm{h}\text{-}\mathrm{proj}(#1)} 
\newcommand{\dgAc}[1]{\mathrm{dgAcy}(#1)}

\newcommand{\scat}[1]{{\mathbf{#1}}} 

\newcommand{\Mod}[1]{#1\text{-}\cat{Mod}} 
\newcommand{\modd}[1]{#1\text{-}\cat{mod}} 
\newcommand{\dgMod}[1]{#1\text{-}\cat{dgMod}} 
\newcommand{\dgCat}{\cat{dgCat}} 

\newcommand{\Hqe}{\Ho{\dgCat}} 

\newcommand{\Spe}{\cat{Spectra}}

\newcommand{\fF}{\fun{F}}

\newcommand\nc {\newcommand}
\newcommand\rnc{\renewcommand}

\nc\sst{\scriptstyle}
\newcommand{\comment}[1]{}
\newcommand{\ri}{\longrightarrow}

\newcommand{\zz}{{\mathbb Z}}

\newcommand{\K}{{\scat K}}
\newcommand{\Ka}{\K^?}
\newcommand{\D}{{\scat D}}
\newcommand{\Ddg}{{\scat D}_{\scat{dg}}}

\nc\op{^{\hbox{\rm\tiny op}}}
\nc\mth{^{\hbox{\rm\tiny th}}}

\nc\script{\mathscr}
\nc\z{\zeta}
\nc\bc{\mathrm{bc}}
\nc\ct{{\script T}}
\nc\cf{{\script F}}
\nc\cg{{\script G}}
\nc\ch{{\script H}}
\nc\ck{{\script K}}
\nc\cl{{\script L}}
\nc\cv{{\script V}}
\nc\ce{{\script E}}
\nc\cs{{\script S}}
\nc\car{{\script R}}
\nc\cd{{\script D}}
\nc\cc{{\script C}}
\nc\ca{{\script A}}
\nc\ci{{\script I}}
\nc\cj{{\script J}}
\nc\co{{\script O}}
\nc\cu{{\script U}}
\nc\cx{{\script X}}
\nc\Cp{{\script P}}
\nc\cq{{\script Q}}
\nc\cy{{\script Y}}
\nc\cz{{\script Z}}
\nc\bd{\begin{description}}
\nc\ed{\end{description}}
\nc\ctob{{\script C}at\big(\ci^{op},\ca\big)}
\nc\clim{{\ds\mathop{\rm lim}_{\ds\longleftarrow}}\,}
\nc\climi{\clim_{\!i}\,}
\nc\climn{\clim^{\!n}\,}
\nc\colim{{\ds\mathop{\rm colim}_{\ds\la}}}
\nc\colimj{{\ds\mathop{\rm colim}_{\ds\la}}{}_{j\,}}
\nc\oa{\overline{\ca}}
\nc\s{\sigma}
\nc\ta{\tau}
\nc\os{\overline\sigma}
\nc\ot{\overline\tau}
\nc\T{\Sigma}
\nc\Tm{\Sigma^{-1}}
\nc\de[1]{{\mathop{\rm deg(#1)}}}
\nc\Ad[1]{\mathop{\rm Ad}(#1)}
\nc\ad[1]{\mathop{\rm ad}(#1)}
\nc\kth{{\it K}--theory}
\nc\loc[1]{{\text{\rm Loc(#1)}}}
\nc\coloc[1]{{\text{\rm Coloc}(#1)}}

\def\der #1 {D\left(#1\right)}
\nc\prf{\begin{proof}}
\nc\eprf{\end{proof}}
\nc\ds{\displaystyle}
\nc\Tor{\text{\rm Tor}}

\nc\cb{{\script B}}
\nc\ab{{\script A}b}

\nc\be{\begin{roenumerate}}
\nc\ee{\end{roenumerate}}

\nc\csab{{\script C}at\big(\cs^{op},\ab\big)}
\nc\ctab{{\script C}at\Big({\{\ct^\alpha\}}^{op},\ab\Big)}
\nc\csex{{\script E}x\big(\cs^{op},\ab\big)}
\nc\ctex{{\script E}x\Big({\{\ct^\alpha\}}^{op},\ab\Big)}
\nc\sub{\qquad\subseteq\qquad}
\nc\ctr[1]{{\left.\ct\left(-,#1\right)\right|}_{\cs}}
\nc\ctrf[2]{{\left.\ct\left(#1,#2\right)\right|}_{\cs}}
\nc\Ctr[1]{{\left.\ct\left(-,#1\right)\right|}_{\ct^\alpha}}
\nc\Ctrf[2]{{\left.\ct\left(#1,#2\right)\right|}_{\ct^\alpha}}

\nc\la{\longrightarrow}
\nc\nin{\noindent}
\nc\cad[1]{\text{card}(#1)}
\nc\eq{\quad=\quad}
\nc\BA{\begin{array}{c}}
\nc\EA{\end{array}}
\nc\barr{
\[
\begin{array}{cccccccccccccccc}
}
\nc\earr{
\end{array}
\]
}
\nc\as[1]{{\langle S\rangle}^{#1}}
\nc\shi{\text{\it shift}}

\nc\yy[1]{{\left.\ct\left(-,#1\right)\right|}_{\ct^c}}
\nc\vrep[2]{{\left.\ct\left(#1,#2\right)\right|}_{\ct^\alpha}}
\nc\da{\downarrow}
\nc\HHom{{\script H}{\mathop{\rm om}}}
\nc\mr\Modtc
\nc\PExt{{\mathop{\rm PExt}}}
\nc\stm{\text{\rm stmod}(kG)}
\nc\stM{\text{\rm StMod}(kG)}
\nc\e{\varepsilon}
\nc\p{\varphi}

\nc\rs{\s^{-1}A}
\nc\br{{\{\s^{-1}A\}}}
\nc\ra\ri
\nc\y[1]{\mathbf{y}#1}
\nc\x[1]{\mathbf{z}#1}
\nc\mmod[1]{#1\text{--\rm mod}}
\nc\MMod[1]{\text{\rm Mod--}#1}
\nc\Md {\ensuremath{\mathop{\textup{Mod}}}}
\rnc\mod[1]{\ensuremath{\mathop{#1\textup{--mod}}}\xspace}
\nc\Modtc{\Mod{\ct^c}}
\nc\pgldim[1]{\mathop{\rm pgldim}\,#1}
\nc\tf{{\rm [TR5]}}
\nc\tfs{{\rm [TR5$^*$]}}
\nc\Fun{\text{\rm Funct}(F\op,\ab)}
\nc\sym{\text{\rm Sym}}
\nc\sgn{\text{\rm sgn}}
\nc\Pro{\text{\rm Prod}^{}_\alpha(F\op,\ab)}
\nc\Yt[1]{{\left.\Hom_\ct^{}\left(-,#1\right)\right|}_F^{}}
\nc\dl{\delta}
\nc\Proj[1]{#1\text{--\rm Proj}}
\nc\proj[1]{#1\text{--\rm proj}}
\nc\Flat[1]{#1\text{--\rm Flat}}
\nc\Inj[1]{#1\text{--\rm Inj}}
\nc\Ima{\mathrm{Im}}
\nc\Ker{\mathrm{Ker}}
\nc\ov{\overline}
\nc\wt{\widetilde}
\nc\wh{\widehat}
\nc\ph{\varphi}
\nc\tstr{{\it t}--structure}
\nc\tstrs{{\it t}--structure }
\nc\spec[1]{{\text{\rm Spec}\left(#1\right)}}

\nc\EProd{\text{\rm EProd}}
\nc\ECoprod{\text{\rm ECoprod}}
\nc\Prod{\text{\rm Prod}}
\nc\ldimp{\text{\rm LDim}^{\prod}}
\nc\ldimc{\text{\rm LDim}^{\coprod}}
\nc\gen[2]{{\langle#1\rangle}^{}_{#2}}
\nc\genu[3]{{\langle#1\rangle}^{[#3]}_{#2}}
\nc\ogen[1]{\ov{\langle#1\rangle}}
\nc\ogenun[2]{\ov{\langle#1\rangle}_{#2}^{}}
\nc\ogenu[3]{\ov{\langle#1\rangle}^{[#3]}_{#2}}
\nc\ogenul[3]{\ov{\langle#1\rangle}^{(-\infty,#3]}_{#2}}
\nc\ogenuf[3]{\ov{\langle#1\rangle}^{[#3,\infty)}_{#2}}
\nc\genuf[3]{{\langle#1\rangle}^{[#3,\infty)}_{#2}}
\nc\genul[3]{{\langle#1\rangle}^{(-\infty,#3]}_{#2}}
\nc\dperf[1]{\D^{\mathrm{perf}}(#1)}
\nc\dcoh{\mathbf{D}^b_{\mathrm{coh}}}
\newcommand{\Dqc}{{\mathbf D_{\text{\bf qc}}}}
\newcommand{\Dqcmi}{{\mathbf D_{\text{\bf qc}}^-}}
\newcommand{\Dqcpl}{{\mathbf D_{\text{\bf qc}}^+}}
\newcommand{\Dqcb}{{\mathbf D_{\text{\bf qc}}^b}}
\newcommand{\Dqcp}{{\mathbf D_{\text{\bf qc}}^p}}
\newcommand{\Dqca}{{\mathbf D_{\text{\bf qc}}^?}}
\newcommand{\Dqcmis}[1]{{\mathbf D_{\text{\bf qc},#1}^-}}
\newcommand{\Dqcpls}[1]{{\mathbf D_{\text{\bf qc},#1}^+}}
\newcommand{\Dqcbs}[1]{{\mathbf D_{\text{\bf qc},#1}^b}}
\newcommand{\Dqcps}[1]{{\mathbf D_{\text{\bf qc},#1}^p}}
\newcommand{\Dqcpbs}[1]{{\mathbf D_{\text{\bf qc},#1}^{p,b}}}
\newcommand{\Dqcpb}{{\mathbf D_{\text{\bf qc}}^{p,b}}}
\newcommand{\Db}{\mathbf{D}^b}
\nc\dmcoh{\mathbf{D}^-_{\mathrm{coh}}}
\nc\dscoh{\mathbf{D}^{}_{\mathrm{coh}}}
\nc\RHHom{{\script{RH}}{\mathrm{om}}}
\nc\Coprod{\mathrm{Coprod}}
\nc\COprod{\mathrm{coprod}}
\nc\add{\mathrm{add}}
\nc\Add{\mathrm{Add}}
\nc\Smr{\mathrm{smd}}
\nc\LL{\mathbf{L}}
\nc\R{\mathbf{R}}
\nc\wi{\wt{\text{\it\i}}}
\nc\exal{\ce\text{\it x}(\ct^\alpha,\ab)}
\nc\exalz{\ce\text{\it x}_{\aleph_0}^{}(\ct^\alpha,\ab)}
\nc\fc{\mathfrak{C}}
\nc\fl{\mathfrak{L}}
\nc\fs{\mathfrak{S}}
\nc\Prf{\text{\bf Perf}}
\nc\Enh{\mathbf{Enh}}
\nc\Tri{\mathbf{Tri}}
\nc\fgt{\mathsf{Fgt}}
\newcommand{\Dqcs}[1]{{\mathbf D}_{\text{\bf qc},#1}}
\newcommand{\Dqcsa}[1]{{\mathbf D^?_{\text{\bf qc},#1}}}
\nc\dcohs[1]{\mathbf{D}^-_{\mathrm{coh},#1}}
\nc\coh{\mathbf{Coh}}
\nc\qc{\mathbf{Qcoh}}
\nc\vect{\mathbf{Vect}}
\nc\dperfs[2]{\D_{#1}^{\mathrm{perf}}(#2)}

\newcommand{\sHom}{\mathcal{H}om} 
\newcommand{\dgHom}{\underline{\sHom}} 

\nc\hoco{
\begin{picture}(40,10)
\put(20,0){\makebox(0,0)[b]{\text{\rm Hocolim}}}
\put(5,-2){\vector(1,0){30}}
\end{picture}\,\,}

\renewcommand{\leq}{\leqslant}
\renewcommand{\geq}{\geqslant}
\nc\tst[1]{\left({#1}^{\leq0},{#1}^{\geq1}\right)}
\nc\tstv[2]{\left({#1}_{#2}^{\leq0},{#1}_{#2}^{\geq1}\right)}
\nc\tsth[2]{{#1}_{#2}^{\heartsuit}}

\nc\holim{
\begin{picture}(40,10)
\put(20,0){\makebox(0,0)[b]{\text{\rm Holim}}}
\put(35,-2){\vector(-1,0){30}}
\end{picture}}

\begin{document}

	\author[A.~Canonaco]{Alberto Canonaco}
\address{A.C.: Dipartimento di Matematica ``F.\ Casorati''\\
        Universit{\`a} degli Studi di Pavia\\
        Via Ferrata 5\\
        27100 Pavia\\
        Italy}
	\email{alberto.canonaco@unipv.it}	
    
  \author[A.~Neeman]{Amnon Neeman}
    \address{A.N.: Dipartimento di Matematica ``F.~Enriques''\\Universit\`a degli Studi di Milano\\Via Cesare Saldini 50\\ 20133 Milano\\ Italy}
    \email{amnon.neeman@unimi.it}
    
    \author[P.~Stellari]{Paolo Stellari}
    \address{P.S.: Dipartimento di Matematica ``F.~Enriques''\\Universit\`a degli Studi di Milano\\Via Cesare Saldini 50\\ 20133 Milano\\ Italy}
    \email{paolo.stellari@unimi.it}
    \urladdr{\url{https://sites.unimi.it/stellari}}

 \title[Weakly approximable categories and enhancements: a survey]
       {Weakly approximable triangulated categories and enhancements: a survey}

 \thanks{A.~C.~is a member of GNSAGA (INdAM) and was partially supported by the research project PRIN 2022 ``Moduli spaces and special varieties''. A.~N.~was partly supported 
   by Australian Research Council Grants DP200102537 and DP210103397,
   and by ERC Advanced Grant 101095900-TriCatApp.
    	P.~S.~was partially supported by the ERC Consolidator Grant ERC-2017-CoG-771507-StabCondEn, by the research project PRIN 2022 ``Moduli spaces and special varieties'', and by the research project FARE 2018 HighCaSt (grant number R18YA3ESPJ)}



\begin{abstract}
This paper surveys some recent results, concerning the intrinsicness of natural subcategories of weakly approximable triangulated categories. We also review the results about uniqueness of enhancements of triangulated categories, with the aim of showing the fruitful interplay. In particular, we show how this leads to a vast generalization of a result by Rickard about derived invariance for schemes and rings. 
\end{abstract}

\maketitle

\setcounter{tocdepth}{1}
\tableofcontents

\section*{Introduction}

If we look at the role played by triangulated categories, in modern algebraic geometry, algebraic topology and representation theory, we realize that they have moved from the periphery to the center. Triangulated categories
started out  being little more than a flexible and general language, to state and frame homological results. But, as of today, they have become one of the key ingredients in proving such results. One example is provided by the theory of stability conditions, which during the last couple of decades has created vast generalizations of the `classical' notion(s) of stability for sheaves.
And this, in turn, has led to amazing applications, ranging from algebraic geometry to string theory and mathematical physics. More recently, the observation that triangulated categories are just the shadow of their higher categorical enhancements, has led to the development of derived algebraic geometry.

We now turn to the question of which developments,
and more specifically which techniques,  have proved most useful
and influential turning triangulated categories from a fringe subject to the
powerful tool they are today. And, among these, three stand out prominently. 

The first is the notion of generation for triangulated categories. Roughly speaking, this is all about the number of steps it takes for one object, or a
specified collection of objects, to generate the entire category. This has led to various notions of dimension for triangulated categories, and was applied to a variety of problems ranging from Brown representability to Grothendieck duality and (conjecturally) to birational geometry.

The second technique is that of \tstr{s} which, among  other things, induce homological functors. Those homological functors take an abstract triangulated category to a simpler abelian category, thus bringing us back to the
realm of classical homological algebra---albeit with new and wonderful coefficients. Going back some paragraphs above, \tstr{s} are one of the key ingredient in the definition of a stability condition.

In this paper we are interested in the interplay between the notions of generation and of \tstr{s}. In a series of papers, this interplay inspired the second author to formalize and explore the notion of \emph{weakly approximable} triangulated categories, as well as the stronger concept of \emph{approximable} triangulated categories. Roughly speaking, a triangulated category is weakly approximable if it possesses a generator which is compact, as well as a \tstr, such that the homology of a bounded above object can be approximated by the generator in a prescribed manner. The formalism of weakly approximable triangulated categories allowed the second author to use the new notion in a number of applications, including the proof of a conjecture about the non-existence of bounded \tstr{s} for geometric categories \cite{Nee2}.

The third technique is that of higher categorical enhancements. Note that, while there is no reason for a weakly approximable triangulated category to have higher categorical enhancements, all those arising from geometric/algebraic/topological contexts do have such enhancements. The question of how canonical or unique such enhancements are has a long history: it was subject to several conjectures, most of which were proved, in great generality, by the three authors in \cite{CNS1}. As explained in \autoref{subsec:enhwip} we expect that the quest, for uniqueness of enhancements for weakly approximable triangulated categories, might be even more successful and provide stronger results than those obtained so far. 

The aim of this survey is to go beyond that, and explain how fruitful the interplay between weakly approximable triangulated categories and higher categorical enhancements can be. In particular we will prove, as an easy application of our results, a vast generalization of a beautiful result by Rickard---see \autoref{thm:Rickard}. The theorem essentially states that, for a ring satisfying mild conditions, all the triangulated categories classically associated to it carry the same amount of algebraic/geometric information. Furthermore, if two such categories are equivalent, then there exists an equivalence of \emph{Morita type.}
This means that we can find an equivalence, possibly different from the one given to us, which is induced by tensor product with a bimodule. We prove a generalization
of Rickard's result in \autoref{thm:main3}. As a special case we have Rickard's old result, but another special case
shows that a parallel statement holds for any quasi-compact and quasi-separated scheme. The take-home message is that, for any quasi-compact and quasi-separated scheme $X$, all the classically associated triangulated categories
\[
\D^?_\mathbf{qc}(X)\qquad\dperf{X}
\]
for $?=b,+,-,\emptyset,p,(p,b)$, carry the same geometric information about $X$. As the notation $\D^p_\mathbf{qc}(X)$ (respectively $\D^{p,b}_\mathbf{qc}(X)$) is not completely standard, let us mention that it coincides with $\D^-_\mathbf{coh}(X)$ (respectively $\D^b_\mathbf{coh}(X)$) when $X$ is noetherian.

The result above is surprising: it goes against the general belief that, when $X$ is not regular and thus the inclusion $\dperf{X}\subseteq \D^{p,b}_\mathbf{qc}(X)$ is proper, then one should carefully choose one of the two categories when dealing with the study of the geometry of $X$. In view of our result this is not the case. Such a choice is a matter of technical convenience, to be determined by practical and computational considerations.

Let us conclude by mentioning that, in \autoref{subsec:appl}, we deduce some interesting applications concerning the category of singularities. It is proved to be a derived invariant. The same is true for the property of being regular, in the case of finite-dimensional, noetherian schemes.

\subsection*{Structure of the paper}

The paper is organized to highlight the idea that the two apparently unrelated theories, of weakly approximable triangulated categories and of uniqueness of enhancements, may fruitfully meet to produce striking new results.

To this end we start out with \autoref{sec:compgen} and \autoref{sec:tstr}, which review the two main ingredients in the definition of weakly approximable triangulated categories: the notion of generation and the one of \tstr. In both section we keep the same structure: we first present the main definitions and properties, and then we focus on the three main examples which play a role in the current manuscript. These are the derived category of modules over a ring, the derived category of complexes of sheaves of $\co_X^{}$--modules, on a scheme $X$,
with quasi-coherent cohomology supported on a closed subset, and the homotopy category of spectra. One special and important construction is the one involving compactly generated \tstr s (see \autoref{subsec:excompgen}).

Following this review of classical material,
in \autoref{sec:wat} we introduce the first main
(relatively recent) key player in this survey: weakly approximable triangulated categories. After recalling the main definitions and properties, we return to our favorite three examples. Finally, in \autoref{subsec:wasubcat}, we come to new work, outlining our main result about the intrinsicness of natural subcategories of weakly approximable triangulated categories---see \autoref{thm:main1}.

In \autoref{sec:enh} we move on to the second key player (somewhat less recent) of this survey: enhancements of triangulated categories and their uniqueness. Again, after reviewing the main definitions and constructions (\autoref{subsec:enhdefex}) we return to examples (\autoref{subsec:dgex}), and then in \autoref{subsec:enhuniq} we present our second main (recent) result, see \autoref{thm:main2}. And then, in \autoref{subsec:enhwip}, we take some time to discuss our work in progress in the direction of generalizing and extending \autoref{thm:main2}. This work in progress constitutes the interplay between uniqueness of
enhancements and the metric methods that are a key to the
study of weakly approximable triangulated categories and their subcategories.

Finally in \autoref{sec:Morita} we combine \autoref{sec:wat} and \autoref{sec:enh} to prove a generalization of a celebrated result by Rickard and some of its geometric applications.

\section{Compactly generated triangulated categories}\label{sec:compgen}

In this section we briefly recall how a triangulated category might be generated by a subcategory, and we focus on the special case when the generating subcategory consists of a single object. Of particular interest is the case where the generating object is compact, see \autoref{subsec:gen}. And then, in \autoref{subsec:excompgen}, we discuss the three classes of examples which we will come back to again and again.

\subsection{How is a category generated?}\label{subsec:gen}

As we have already mentioned, in this paper we are mainly interested in triangulated categories which are suitably generated. In this  section we set the notation which will be used throughout the paper.

We start with a triangulated category $\ct$ and two full subcategories $\ca$ and $\cb$. There are three natural operations we can perform: taking extensions, taking coproducts and taking direct summands. These operations will be encoded as follows:
\begin{itemize}
\item
$\ca*\cb\subseteq\ct$ is the full subcategory of all objects $C\in\ct$ for which there exists a distinguished triangle $A\la C\la B$ with $A\in\ca$ and $B\in\cb$.
\item
$\add(\ca)\subseteq\ct$ is the full subcategory whose objects are all finite direct sums of objects of $\ca$. In case $\ct$ has (set indexed) coproducts, then $\Add(\ca)\subseteq\ct$ is the full subcategory whose objects are all the small coproducts of objects of $\ca$. 
\item
$\Smr(\ca)\subseteq\ct$ is the full
subcategory with objects all direct summands of objects of $\ca$.
\setcounter{enumiv}{\value{enumi}}
\end{itemize}
Using these operations we can define the full subcategory $\COprod(\ca)$. It is the smallest full subcategory $\cs\subseteq\ct$ satisfying
\[
\ca\subseteq\cs,\qquad\cs*\cs\subseteq\cs,\qquad\add(\cs)\subseteq\cs.
\]
And we also set  
\[
\gen\ca{}:=\Smr\left(\COprod(\ca)\right).
\]

If we are in the situation where $\ct$ has coproducts, the parallel definitions set $\Coprod(\ca)$ to be the
smallest full subcategory $\cs\subseteq\ct$ satisfying
\[
\ca\subseteq\cs,\qquad\cs*\cs\subseteq\cs,\qquad\Add(\cs)\subseteq\cs.
\]
And, again in parallel, we set
\[
\ogen\ca{}:=\Smr\left(\Coprod(\ca)\right).
\]

We will be interested in special full subcategories $\ca$, namely those containing a single object $G$ and some specified shifts of it. In symbols, if $m\leq n$ are integers, possibly infinite, then we set
\[
G[m,n]:=\{\sh[i]{G}\st i\in\zz\text{ and } m\leq -i\leq n\}.
\]
And then we will use the shorthand
\[
\genu G{}{m,n}:=\gen{G[m,n]}.
\]

\begin{remark}\label{rmk:genfin1}
In the special case $m=-\infty$ and $n=+\infty$, it is easy to see that we get the equality
\[
\genu G{}{-\infty,+\infty}=\bigcup_{i>0}\genu G{}{-i,i}.
\]
\end{remark}

Of course, we get analogue constructions when $\ct$ has coproducts:
\[
\ogenu G{}{m,n}:=\ogen{G[m,n]}.
\]

Let us now stick to the case of a triangulated category $\ct$ with coproducts. An object $C\in\ct$ is \emph{compact} if, given any collection $\{C_i\}_{i\in I}$ of objects in $\ct$, the natural morphism
\begin{equation}\label{eqn:cmpt}
\coprod_{i\in I}\Hom(C,C_i)\lto\Hom\left(C,\coprod_{i\in I}C_i\right)
\end{equation}
is an isomorphism\footnote{When dealing with stable $\infty$-categories one usually requires that an object is compact if it satisfies a similar property where coproducts are replaced by filtered colimits. For triangulated categories, filtered colimits should be replaced by homotopy colimits. It is a classical result (see, for example, \cite[Lemma 3.12]{R}) that a compact object in the sense of the present paper is such that the analogue of \eqref{eqn:cmpt} with homotopy colimits in place of coproducts is an isomorphism.}. We denote by $\ct^c$ the full subcategory of $\ct$ consisting of compact objects.

\begin{remark}\label{rmk:comp1}
It is not hard to see, using the definition, that $\ct^c$ is always triangulated. But there is no reason for $\ct^c$ to be in any sense ``big'': it might very well be trivial, see \cite{K1} for a simple example. We will go back later to additional `pathologies' of the subcategory $\ct^c$.
\end{remark}

\begin{definition}\label{def:compgen}
A triangulated category $\ct$ with coproducts is \emph{compactly generated} by $G\in\ct^c$ if $\ct={\ogen G{}}^{[-\infty,+\infty]}$.
\end{definition}

\begin{remark}\label{rmk:equivcomp}
  For $\ct$ with coproducts one can give an alternative definition of compact generation. It turns out to be equivalent to require that, if $E\in\ct$ is such that $\Hom(G[-\infty,+\infty],E)=0$, then $E\cong 0$. It is easy to see that the condition in \autoref{def:compgen} implies the one above. Going in the
  other direction is an application of Brown representability, which holds for $\ct$ satisfying the above vanishing condition. This can be found in several
  places in the literature (see, for example, \cite[Remark 2.2]{CS:surv1}).
\end{remark}

Note that, by \cite{Nee1}, if $\ct$ has a compact generator $G$, meaning $\ct={\ogen G{}}^{[-\infty,+\infty]}$, then $\ct^c=\genu G{}{-\infty,+\infty}$. When this happens, we say that $G$ \emph{classically generates} $\ct^c$.

\begin{remark}\label{rmk:genfin2}
Assume now that $G_1$ and $G_2$ are two compact generators for $\ct$. By \autoref{rmk:genfin1}, there exists a positive integer $M$ such that both equalities
\[
G_1\in\genu {G_2}{}{-M,M}\qquad G_2\in\genu {G_1}{}{-M,M}
\]
hold.
\end{remark}

\subsection{Examples}\label{subsec:excompgen}

In this presentation, we are mainly interested in three sets of examples. Let us discuss them in some detail.

\subsubsection*{Modules over a ring}
 Let $R$ be a (possibly noncommutative) ring. We will consider the Grothendieck abelian category $\Mod{R}$ of left $R$-modules and its unbounded derived category $\D(\Mod{R})$, which is a triangulated category with coproducts. Since $\Hom_{\D(\Mod{R})}(R,-)\cong H^0(-)$, it is clear that the object $R$ is compact in $\D(\Mod{R})$. Furthermore, given any complex $D\in\D(\Mod{R})$, if $\Hom(\sh[i]{R},D)=0$ for all $i\in\ZZ$, then $D\cong 0$. Thus, by \autoref{rmk:equivcomp}, $R$ is a compact generator for $\D(\Mod{R})$. Following the standard notation, we set $\dperf{R}:=\D(\Mod{R})^c$. As $\dperf{R}=\genu R{}{-\infty,+\infty}$, one can easily get the alternative description
 \[
 \dperf{R}=\K^b(\proj R),
 \]
 that is the homotopy category of bounded complexes of finitely generated projective $R$-modules. Observe that the compact generator $R$ satisfies the strong vanishing $\Hom(R,\sh[i]{R})=0$ for $i>0$.

 Assume now that $R$ is \emph{left coherent}---this means that every finitely generated left ideal is finitely presented. In this case the category $\modd{R}$ of finitely presented $R$-modules is abelian, and we have a natural identification
 \[
 \Db(\modd{R})=\K^{-,b}(\proj R),
 \]
 where the category on the right-hand side is the homotopy category of bounded above complexes of finitely generated projective $R$-modules with bounded cohomology. In particular, we get inclusions
 \[
 \dperf{R}\subseteq\Db(\modd{R})\subseteq\D(\Mod{R}).
 \]
Note that any left noetherian ring is left coherent (but not conversely). 

\subsubsection*{Schemes with support}

Let us now move to the purely commutative and geometric side and assume that $X$ is a quasi-compact and quasi-separated scheme. Let $Z\subseteq X$ be a closed subscheme such that $X\setminus Z$ is quasi-compact. Consider the triangulated category $\Dqcs Z(X)$ consisting of unbounded complexes of
sheaves of $\mathcal{O}_X$-modules with quasi-coherent cohomology sheaves supported on $Z$.

This triangulated category has coproducts and, by a beautiful result by Rouquier \cite[Theorem 6.8]{R}, it is generated by an object $G$ that is not only compact in $\Dqcs Z(X)$ but is actually contained in $\Dqcs Z(X)\cap\Dqc(X)^c$. As a consequence, we have natural identifications
\[
\Dqcs Z(X)^c=\Dqcs Z(X)\cap\Dqc(X)^c=\Dqcs Z(X)\cap\dperf{X}=\dperfs ZX.
\]
Here we use the fact that, for a quasi-compact and quasi-separated scheme $\Dqc(X)^c=\dperf{X}$, where the latter category consists of complexes which are locally quasi-isomorphic to bounded complexes of locally free sheaves of finite rank. Clearly $\dperfs ZX$ stands for the full subcategory of $\dperf{X}$ consisting of complexes with cohomology supported on $Z$.

\begin{remark}\label{rmk:ideaproof1}
The idea of the proof is relatively easy to explain. Since $X$ is quasi-compact and quasi-separated, one can reduce the problem to $X$ affine and $Z$ described by finitely many equations $f_1=\dots=f_n=0$. One then proves that the complex
\[
\bigotimes_{i=1,\dots n}\left(0\lto\mathcal{O}_X\stackrel{f_i}{\lto}\mathcal{O}_X\lto 0\right)
\]
is a generator of $\Dqcs Z(X)$, which is clearly compact in $\Dqc(X)$. Note that this result has been reproved more recently in \cite{Nee2}.
\end{remark}

Finally, as the generator $G$ belongs to
$\dperf{X}$, the general result \cite[Tag 09M4]{Stack} implies that
\begin{equation}\label{eq:vanHomG}
\Hom(G,\sh[i]{G})=0,
\end{equation}
for $i\gg 0$.

\subsubsection*{Homotopy category of spectra}

Next we discuss another triangulated category, of a more topological flavor, for which the theory discussed in the next few sections applies. First recall:
for a topological space $X$ with a basepoint $x_0\in X$, the
pointed space $\Sigma X$ is obtained from $X\times [0,1]$ by collapsing
\[
X\times\{0\}\bigcup X\times\{1\}\bigcup\Big(\{x_0\}\times[0,1]\Big)
\]
to a single point, which is declared to be the basepoint of $\Sigma X$. The space $\Sigma X$
is called the \emph{reduced suspension} of $X$.
The reduced suspension defines an endofunctor $\Sigma$ on the category of pointed topological spaces, and it is known to have a right adjoint $\Omega$. This right adjoint takes a pointed space $X$ to the pointed space $\Omega X$, of pointed loops on $X$.

Next we come to the category
we want to consider, the category of spectra. 
It can be constructed as follows: an object is a
so-called
\emph{$\Omega$--spectrum.} An $\Omega$--spectrum is a sequence $F$ of pointed topological spaces $F_n$, together with pointed continuous maps $\sigma_n\colon\Sigma F_n\to F_{n+1}$. The continuous maps
$\sigma_n\colon\Sigma F_n\to F_{n+1}$ are called the \emph{structural morphisms} of the spectrum. The
structural morphisms must satisfy the following condition: 
the map corresponding to $\sigma_n:\Sigma F_n\la F_{n+1}$,
under the adjunction $\Sigma\dashv\Omega$, has to be a
weak homotopy equivalence
$F_n\la\Omega F_{n+1}$. In case the
reader is curious: the weak homotopy equivalences
$F_n\la\Omega F_{n+1}$ is the reason for calling this
incarnation an $\Omega$--spectrum; the homotopy
category of spectra, which is what we are constructing,
has other, equivalent descriptions where the objects
are spectra satisfying less restrictive hypotheses. Anyway, 
a morphism of $\Omega$--spectra $E\to F$ is a sequence of pointed continuous maps $E_n\to F_n$, strictly compatible with the structural morphisms. This definition yields $\Spe$, \emph{the category of $\Omega$--spectra}.

We can then define the homotopy groups of a spectrum $F$ as follows
\[
\pi_n(E):=\colimm{i\geq\max\{0,-n\}}\pi_{n+i}(E_i).
\]
Here the colimit is with respect to the natural isomorphisms
\[
\pi_{n+i}(E_i)\iso\pi_{n+i}(\Omega E_{i+1})\iso\pi_{n+i+1}(E_{i+1}).
\]
A \emph{weak equivalence} is a morphism between $\Omega$--spectra that induces an isomorphism on the homotopy groups. The \emph{(stable) homotopy category of spectra} $\ho{\Spe}$ is defined as the localization of $\Spe$ by the class of the weak equivalences. The category $\ho{\Spe}$ is a $\ZZ$-linear triangulated category with coproducts.

As we have already said this category has other descriptions,
in which the objects are different. Thus, for example, there
is no need to assume that the maps
$\sigma_n:\Sigma F_n\la F_{n+1}$
correspond by adjunction to weak homotopy equivalences
$F_n\la\Omega F_{n+1}$. We can start with
any spectrum, meaning any sequence $F$ of pointed spaces $F_n$,
and place no restrictions on the structural morphisms
$\sigma_n:\Sigma F_n\la F_{n+1}$. And then
we turn it into an $\Omega$--spectrum by setting
\[E_n=\colimm{i\to\infty}\Omega^{i-n} F_i\ .\]
In this second description, the following
is an important example of an object. The \emph{sphere spectrum} $S^0$, which is \emph{not} an
$\Omega$--spectrum, is obtained by taking $F_n$ to be the $n$-dimensional (pointed) sphere $S^n$, and
the structural morphism $\Sigma F_n\la F_{n+1}$
is the homeomorphism $\Sigma S^n\to S^{n+1}$. It turns out that the object $S^0$ is a compact generator of $\ho{\Spe}$. There is a very brief discussion in \cite[Example 3.2]{Nee3}, and much fuller accounts in the homotopy theoretic literature.

\section{\tstr s}\label{sec:tstr}

In this section we introduce the second ingredient in the definition of weakly approximable triangulated categories: \tstr s. After a brief discussion of the definition and the basic properties in \autoref{subsec:deftstr}, in \autoref{subsec:comptstr} we discuss the crucial construction of compactly generated \tstr s, and finally we return to the main examples in \autoref{subsec:extstr}.

\subsection{Definition and basic properties}\label{subsec:deftstr}

Let us start with the main definition of this section. Our presentation is terse, but we refer to Sections 2.1 and 2.2 in \cite{CNS2} for a detailed discussion and comparison between different approaches and descriptions. For later use, recall that a subcategory $\cc$ of
a category $\cd$ is \emph{strictly full} if it is full and,
given an object
$A\in\cc$ and an isomorphism $A\la B$ in $\cd$, then $B$ is in $\cc$.

\begin{definition}\label{def:tstr}
Let $\ct$ be a triangulated category. A \emph{\tstr}\ on $\ct$ is a pair
of strictly full subcategories $\tau=(\ct^{\leq 0},\ct^{\geq 1})$ satisfying:
\begin{itemize}
\item[(i)]
$\sh{\ct^{\leq 0}}\subseteq\ct^{\leq 0}$ and $\sh[-1]{\ct^{\geq 1}}\subseteq\ct^{\geq 1}$.
\item[(ii)]
$\Hom(\ct^{\leq 0},\ct^{\geq 1})=0$.
\item[(iii)]
For any object $B\in\ct$ there exists a distinguished triangle $A\la B\la C\la\sh{A}$, 
with $A\in\ct^{\leq 0}$ and $C\in\ct^{\geq 1}$.  
\end{itemize}
\end{definition}

In presence of a \tstrs we can, for all $n\in\ZZ$,
define the following
auxiliary strictly full subcategories
\[
\ct^{\leq n}:=\sh[-n]{\ct^{\leq 0}}\ ,\qquad\ct^{\geq n}:=\sh[-n+1]{\ct^{\geq 1}}\ .
\]

\begin{remark}\label{rmk:tstradjoints}
It is clear from the definition that, more generally, given an object $B\in\ct$ and $n\in\ZZ$, there is a distinguished triangle $A\la B\la C\la\sh{A}$, 
with $A\in\ct^{\leq n}$ and $C\in\ct^{\geq n+1}$. It is also not hard to prove that such a triangle is uniquely determined, up to isomorphism. It follows that there are functors
\[
\tau^{\leq n}\colon\ct\lto\ct^{\leq n}\qquad\tau^{\geq n}\colon\ct\lto\ct^{\geq n},
\]
for all $n\in\ZZ$, such that $\tau^{\leq n}(B)=A$ and $\tau^{\geq n+1}(B)=C$. Furthermore, $\tau^{\leq n}$ is the right adjoint to the inclusion $\ct^{\leq n}\hookrightarrow\ct$, while $\tau^{\geq n}$ is the left adjoint to the inclusion $\ct^{\geq n}\hookrightarrow\ct$.
\end{remark}

\begin{ex}\label{ex:degtstr}
It is relatively easy to construct \tstr s on a triangulated category $\ct$. One can simply set $\ct^{\leq 0}=\ct$ and $\ct^{\geq 1}=\{0\}$. Clearly $\ct=\bigcap_n\ct^{\leq n}$, which is non-trivial when $\ct$ is not.
\end{ex}

The example above is an instance of the so called degenerate \tstr s, which will not be relevant in this presentation. Recall that a \tstrs $\tau$ is \emph{non-degenerate} if the intersections of all the $\ct^{\leq n}$ and of all the $\ct^{\geq n}$ are trivial. The aim of \autoref{subsec:comptstr} and \autoref{subsec:extstr} is to provide interesting examples of non-degenerate \tstr s.

A remarkable property which follows from the functorial properties discussed in \autoref{rmk:tstradjoints} is that, for any $n\in\ZZ$, the subcategory $\ct^{\leq n}$ is closed under direct summands and the coproducts which are defined in $\ct$, while $\ct^{\geq n}$ is closed under direct summands and the products that are defined in $\ct$.

A triangulated category with a \tstrs always contains an abelian category. Indeed, the strictly full subcategory 
\[
\ct^\heartsuit:=\ct^{\leq 0}\cap\ct^{\geq 0}
\]
turns out to be an abelian category. It is usually called the \emph{heart} of the \tstr. Note that we get a functor
\begin{equation}\label{eq:cohtstr}
\ch^0:=\tau^{\geq0}\circ\tau^{\leq 0}\colon\ct\la\tsth\ct{},
\end{equation}
which is cohomological, meaning that it takes distinguished triangles in $\ct$ to exact sequences in $\ct^\heartsuit$. One can also define $\ch^n(-):=\ch^0(\sh[n]{-})$, for all $n\in\ZZ$. Given $A\in\ct$, we should think of $\ch^n(A)$ as the $n$-th cohomology of $A$ with respect to the \tstr.

We will later see examples that a single triangulated category may have more than one, distinct \tstr s. It becomes important to be able to compare them.

\begin{definition}\label{def:tstreq}
Two \tstr s $\tau_1=(\ct_1^{\leq 0},\ct_1^{\geq 1})$ and $\tau_2=(\ct_2^{\leq 0},\ct_2^{\geq 1})$ on a triangulated category $\ct$ are \emph{equivalent} if there exists an integer $N\geq 0$ such that
\[
\ct_2^{\leq -N}\subseteq\ct_1^{\leq 0}\subseteq \ct_2^{\leq N}.
\]
\end{definition}

Some comments are in order here.

\begin{remark}\label{rmk:tstr}
(i) First of all, using the definition of \tstr, it is easy to check that the roles of $\tau_1$ and $\tau_2$ in \autoref{def:tstreq} are interchangeable. Furthermore, as suggested by the definition, the relation described above is an equivalence relation between \tstr s. Thus we can look at equivalence classes of \tstr s on a given triangulated category.

(ii) \autoref{def:tstreq} can be equivalently reformulated in terms of the full subcategories $\ct_1^{\geq 1}$ and $\ct_2^{\geq 1}$ and their shifts.

(iii) If $\tau_1=(\ct^{\leq 0}_1,\ct^{\geq 1}_1)$ is a \tstrs on a triangulated category $\ct$ and $\fF\colon\ct\la\ct'$ is a triangulated equivalence, then we get an induced \tstrs $\fF(\tau_1)=(\fF(\ct^{\leq 0}_1),\fF(\ct^{\geq 1}_1))$ on $\ct'$. Clearly, if $\tau_2$ is a second \tstrs on $\ct$ in the same equivalence class as $\tau_1$, then $\fF(\tau_1)$ and $\fF(\tau_2)$ are in the same equivalence class.
\end{remark}

\subsection{Compactly generated \tstr s}\label{subsec:comptstr}

In this section we want to describe an efficient and general procedure to produce \tstr s. To this end let us assume, for the rest of this section, that $\ct$ is a triangulated category with coproducts. Let $\ca\subseteq\ct^c$ be an essentially small full subcategory. We can then form the pair of full subcategories
\[
\tau_\ca:=\big(\Coprod(\ca),\Coprod(\ca)^\perp\big).
\]
Recall that, given a full subcategory $\cs\subseteq\ct$, the collection of objects of the full subcategory $\cs^\perp\subseteq\ct$ is $\{A\in\ct\st\Hom(S,A)=0\text{ for all }S\in\cs\}$. We can now state the following.

\begin{thm}\label{thm:compgentstr}
In the situation above and under the additional assumption that $\sh{\ca}\subseteq\ca$, the pair $\tau_\ca$ is a \tstrs on $\ct$.
\end{thm}

This result originated in \cite[Theorem~A.1 and Proposition~A.2]{Alonso-Jeremias-Souto03}. A different proof was provided by \cite[Theorem 2.3.3]{CNS2}. The reason we felt the need to reprove the result is explained in \cite[Remark 2.3.4]{CNS2}, and boils down to the fact that the original proof in \cite{Alonso-Jeremias-Souto03} seems to implicitly assume that the triangulated category $\ct$ has a higher categorical model (see \autoref{def:dgpreptr} for a precise definition of such a model).

\begin{definition}\label{def:compgentstr}
The \tstrs of the form
$\tau_\ca:=\big(\Coprod(\ca),\Coprod(\ca)^\perp\big)$, where
$\ca\subseteq\ct^c$ is an essentially small full subcategory
satisfying $\sh\ca\subseteq \ca$, will be called the
\emph{\tstrs on $\ct$ compactly generated by $\ca$}.    
\end{definition}

\begin{remark}\label{rmk:preseq}
Following \autoref{rmk:tstr}(iii) note that, since any triangulated equivalence $\fF\colon\ct_1\la\ct_2$ restricts to an equivalence between $\ct_1^c$ and $\ct_2^c$, if $\tau_\ca$ is a \tstrs on $\ct_1$ which is compactly generated by $\ca\subseteq\ct_1^c$, then $\fF(\tau_\ca)$ is a \tstrs on $\ct_2$ which is compactly generated by $\fF(\ca)\subseteq\ct_2^c$.
\end{remark}

Given a \tstrs on $\ct$ which is compactly generated by $\ca$, we can define
\[
\ct^{\leq0}_{\ca,c}:=\ct^c\cap\ct^{\leq0}_\ca \qquad\tsth\ct{\ca,c}:=\ch_\ca^{0}\big(\ct^{\leq0}_{\ca,c}\big).
\]
Here $\ct^{\leq0}_\ca:=\Coprod(\ca)$ and $\ch_\ca^{0}$ denotes the functor \eqref{eq:cohtstr} associated to the \tstrs $\tau_\ca$.

The key properties of $\tau_\ca$ are summarized in the following result, which appeared in \cite[Theorem~8.31]{Saorin-Stovicek20}. See also \cite[3.0.1]{CNS2} for a different proof.

\begin{thm}\label{thm:compgenGroth}
Let $\ct$ be a triangulated category with coproducts,
and let $\ca\subseteq\ct^c$ be an essentially small full subcategory satisfying $\sh{\ca}\subseteq \ca$.
Then the heart $\tsth\ct A$ of the \tstr\ $\tau_\ca=\big(\Coprod(\ca),\Coprod(\ca)^\perp\big)$ is a locally
finitely presented Grothendieck abelian category.
Moreover, $\tsth\ct{A,c}\subseteq\tsth\ct A$
is precisely the full subcategory of
finitely presented objects in $\tsth\ct A$.
\end{thm}

For the non-expert reader, we recall that an abelian category $\cb$ is a \emph{Grothendieck category} if
\begin{itemize}
\item It has small coproducts;
\item Filtered colimits of exact sequences are exact in $\cb$;
\item It has a set of generators, meaning a set $\cs$ of objects in $\cb$ such that, for any $C\in\cb$ there is an epimorphism $S\epi C$, where $S$ is a coproduct of objects in $\cs$.
\end{itemize}
An object $B\in\cb$ is \emph{finitely presented} if the functor $\Hom_\cb(B,-)$ commutes with filtered colimits, and $\cb$ is \emph{locally finitely presented} if the full subcategory of finitely presented objects in $\cb$ is essentially small, and all objects of $\cb$ are filtered colimits of finitely presented objects.

Next we introduce the cases that will be
important to us, where $\ca$ is of a special type. Assume that $\ct$ is compactly generated, and let $G$ be a compact generator. Consider the subcategory $\ca:=G[-\infty,0]$ and set $\tau_G:=\tau_\ca$. For this special $\ca$ we will use the shorthand $\ct_G^{\leq 0}=\ct_\ca^{\leq 0}$.

Now let $H$ be another compact generator for $\ct$. By \autoref{rmk:genfin2}, coupled with the explicit description of $\ct_G^{\leq 0}$ and $\ct_H^{\leq 0}$ in \autoref{def:compgentstr}, it is not hard to see that $\tau_G$ and $\tau_H$ are in the same equivalence class of \tstr s on $\ct$. Thus the following definition makes sense.

\begin{definition}\label{def:prefeqclass}
Let $\ct$ be a triangulated category with coproducts and generated by a single compact object. The equivalence class of \tstr s on $\ct$ which contains $\tau_G$, where $G$ is any compact generator, is the \emph{preferred equivalence class} of \tstr s on $\ct$.
\end{definition}

\begin{remark}\label{rmk:preststr2}
Assume that $\ct_1$ is a triangulated category with coproducts and generated by a single compact object. By \autoref{rmk:tstr}(iii) and \autoref{rmk:preseq}, any triangulated equivalence $\fF\colon\ct_1\la\ct_2$ preserves the preferred equivalence class of \tstr s.
\end{remark}

Next we investigate if the the preferred equivalence class of \tstr s contains examples we recognize. This will be the subject of the next section.

\subsection{Examples}\label{subsec:extstr}

Let us now explore further the three examples in \autoref{subsec:excompgen}: the derived category of modules over a ring, the derived category of complexes of sheaves of $\co_X$--modules over a scheme with quasi-coherent cohomology supported on a closed subset and the homotopy category of spectra.

By the discussion in the previous section, and the observation from \autoref{subsec:excompgen} that $\D(\Mod{R})$, $\Dqcs Z(X)$ and $\ho{\Spe}$ are compactly generated by a single object, we know that these three classes of triangulated categories possess a well-defined preferred equivalence class of \tstr s.

But when $R$ is a ring $\D(\Mod{R})$ is endowed with the standard \tstrs $\tau_\text{stan}$. It is a pair of full subcategories, whose first component is the subcategory $\D^{\leq 0}(\Mod{R})$ consisting of complexes with cohomology in non-positive degrees. Now, since $R$ is a compact generator for $\D(\Mod{R})$, by the discussion in the previous section we can also consider $\tau_R=(\D(\Mod{R})_R^{\leq 0},\D(\Mod{R})_R^{\geq 1})$. By the explicit description of $\D(\Mod{R})_R^{\leq 0}$,
in the paragraphs right before
\autoref{def:prefeqclass}, it is clear that we have the inclusion
\[
\D(\Mod{R})_R^{\leq 0}\subseteq\D^{\leq 0}(\Mod{R}).
\]
On the other hand, every complex $E\in\D^{\leq 0}(\Mod{R})$ is quasi-isomorphic to an (unbounded below) complex
\[
\dots R^{\oplus J_i}\la R^{\oplus J_{i-1}}\la\dots\la R^{\oplus J_1}\la R^{\oplus J_0}\la 0,
\]
where $J_k$ is a set and $R^{\oplus J_k}$ sits in degree $-k$. Such a complex is quasi-isomorphic to the homotopy colimit of its stupid truncations which live in $\Coprod(R[-\infty,0])$. But we observed that $\Coprod(R[-\infty,0])$ is closed under coproducts, and hence under homotopy colimits---these exist in $\D(\Mod{R})$. Thus $E\in\Coprod(R[-\infty,0])$ and we can conclude that
\begin{equation}\label{eq:eq1}
\D(\Mod{R})_R^{\leq 0}=\D^{\leq 0}(\Mod{R}).
\end{equation}
This implies that $\tau_\text{stan}$ is in the preferred equivalence class of \tstr s.

Let us now consider the case of $\Dqcs Z(X)$, where $X$ is a quasi-compact and quasi-separated scheme, containing a closed subscheme $Z$ whose complement $X\setminus Z$ is quasi-compact. Any such category carries a standard \tstrs $\tau_\text{stan}$, whose first component is the subcategory $\Dqcs{Z}^{\leq 0}(X)$, consisting of complexes with cohomology sheaves only in non-positive degrees. The existence of such a special \tstrs leads to the definition of classical and useful full triangulated subcategories of $\Dqcs Z(X)$ (other than $\dperfs{Z}{X}=\Dqcs Z(X)^c)$:
\begin{itemize}
\item
$\Dqcmis Z(X)\subseteq\Dqcs Z(X)$: the
full subcategory of all
complexes $C$ with $\ch^n(C)=0$ for $n\gg0$.
\item
$\Dqcpls Z(X)\subseteq\Dqcs Z(X)$: the
full subcategory of all
complexes $C$ with $\ch^n(C)=0$ for $n\ll0$.
\item
$\Dqcbs Z(X)\subseteq\Dqcs Z(X)$: the
full subcategory of all
complexes $C$ with $\ch^n(C)=0$ for
all but finitely many $n\in\zz$.
\item
$\Dqcps Z(X)\subseteq\Dqcs Z(X)$: the
full subcategory of all
\emph{pseudocoherent} complexes supported on $Z$.
A complex is pseudocoherent if it is locally
isomorphic to a bounded-above complex of
finite-rank vector
bundles.
\item
$\Dqcpbs Z(X)\subseteq\Dqcps Z(X)$: the
full subcategory of all
objects $C\in\Dqcps Z(X)$
such that $\ch^n(C)=0$ for all but finitely
many $n\in\zz$.  
\end{itemize}

Of course, we still have to deal with the question whether $\tau_\text{stan}$ is in the preferred equivalence class. The easy check above involving modules over a ring becomes now a highly non-trivial result.

\begin{thm}[\cite{Nee2}, Theorem 3.2 (iii)]\label{thm:supstantstr}
If $X$ and $Z$ are as above, the \tstrs $\tau_\mathrm{stan}$ is in the preferred equivalence class.
\end{thm}

One easy but important consequence of this result is that all the triangulated subcategories of $\Dqcs Z(X)$ listed above could be equivalently defined by using any of the \tstr s in the preferred equivalence class.

\begin{remark}\label{rmk:vanHoms}
It is clear that the compact generator $R$ of $\D(\Mod{R})$ satisfies the vanishing
\[
\Hom(R,\D(\Mod{R})_R^{\leq -1})=\Hom(R,\D^{\leq -1}(\Mod{R}))=0.
\]
On the other hand, if $G$ is the compact generator for $\Dqcs Z(X)$ mentioned in \autoref{subsec:excompgen}, then there is an integer $n$ such that
\[
\Hom(G,\Dqcs Z(X)_G^{\leq n})=\Hom(G,\Dqcs{Z}^{\leq n}(X))=0.
\]
The vanishing of the first $\Hom$-space follows from \eqref{eq:vanHomG} and the definition of $\Dqcs Z(X)_G^{\leq n}$ while the vanishing of the second $\Hom$-space follows directly from \autoref{thm:supstantstr}.
\end{remark}

Finally, $\ho{\Spe}$ has a standard \tstrs $\tau_\text{stan}$ as well. In this case, the first subcategory in $\tau_\text{stan}$ consists of the so called \emph{connective spectra}. It consists of those spectra $F$ such that $\pi_i(F)=0$ for all $i<0$. As explained in \cite[Example 3.2]{Nee2}, $\tau_\text{stan}$ is in the preferred equivalence class.

Putting it all together, we obtain the following result.

\begin{cor}\label{cor:standtstr}
The standard \tstr s on $\D(\Mod{R})$, $\Dqcs Z(X)$ and $\ho{\Spe}$, for $R$, $X$ and $Z$ as above, are all in the preferred equivalence class.
\end{cor}

\section{Weakly approximable triangulated categories}\label{sec:wat}

This section deals with the notion of weakly approximable triangulated categories, and aims to explain \autoref{thm:main1}, the first main result of this survey paper. Several beautiful properties of these triangulated categories are discussed, with special attention to the three main examples from \autoref{subsec:excompgen}.

\subsection{Definition and basic properties}\label{subsec:wadef}

The crucial definition in this paper is the following, which was first given in \cite[Definition~0.21]{Nee3}. The definition combines in an intriguing way the notion of compact generation and the notion of \tstr. The interplay between the two will be investigated more fully in this section.

\begin{definition}\label{def:wa}
A triangulated category with coproducts $\ct$ is
\emph{weakly approximable} if there exist an integer $A>0$,
a compact generator $G\in\ct^c$ and a \tstrs
$\tau=\tst\ct$ on $\ct$, such that
\begin{itemize}
\item[(i)]
$G\in\ct^{\leq A}$ and $\Hom(G,\ct^{\leq-A})=0$.
\item[(ii)]
For any object $F\in\ct^{\leq0}$ there exists in $\ct$ a distinguished triangle
$E\la F\la D\la\sh{E}$ with $E\in\ogenu G{}{-A,A}$ and with
$D\in\ct^{\leq-1}$.
\end{itemize}
\end{definition}

For the sake of conciseness, a triple $(G,\tau,A)$ satisfying (i) and (ii) in \autoref{def:wa} will be referred to as \emph{weak approximation data} for $\ct$.

\begin{remark}\label{rmk:approx}
Following \cite{Nee3}, one can also define \emph{approximable} triangulated categories. In this variant one replaces, in (ii) above, the condition $E\in\ogenu G{}{-A,A}$ with the more restrictive one $E\in\ogenu G{A}{-A,A}$. We did not introduce the symbol $\ogenu G{A}{-A,A}$ in \autoref{subsec:gen}, as it will not be relevant for the rest of the current survey. Roughly, the subscript $A$ prescribes that the object $E$ can be obtained by at most $A$ iterations of the three basic operations involved in the definition of $\ogenu G{}{-A,A}$. In the next section we will point out that there are examples of triangulated categories which are weakly approximable but not approximable.
\end{remark}

Axiom (ii) in the definition above is crucial to many of the proofs and motivates the terminology. Given a bounded above complex in a weakly approximable triangulated category, then (ii) says that its top cohomology can be approximated by an object generated using uniformly bounded shifts of the generator.

That said, it is an open question whether (ii) is
implied by (i). The evidence that this may be the case is a result
due to Bondarko and Vostokov:

\begin{remark}\label{rmk:wadep}
Assume that $\ct$ is triangulated category with coproducts and a single compact generator $G$. Assume further that
\begin{equation}\label{eq:vancomp}
\Hom(G,G[i])=0
\end{equation}
for $i\geq 2$. We will see in the next section that, although this is a restrictive assumption, it is satisfied in several interesting situations.

Next we sketch the proof that, for such a triangulated category, (i) in \autoref{def:wa} implies (ii) in the same definition. Take the \tstrs $\tau_G$. It is clear that the triple $(G,\tau_G,2)$ satisfies \autoref{def:wa}(i). Let us now show that, assuming \eqref{eq:vancomp}, such a triple yields weak approximation data for $\ct$. This means that, given $F\in\ct_G^{\leq 0}$, we want to find a distinguished triangle
\[
E\la F\la D,
\]
where $D\in\ct_G^{\leq -1}$ and $E\in\ogenu G{}{-2,2}$. In fact we will prove a stronger statement: we will show that
$E$ may be chosen to lie in $\ogenu G{}{0,0}=:\ct^0$. 

Let us denote by $\cs$ the collection of objects $F\in\ct_G^{\leq 0}$ for which such a triangle exists with $E\in\ct^0$. Clearly, $\cs$ contains all the shifts $\sh[i]{G}$, for $i\geq 0$. Furthermore, it is closed under coproducts. In order to prove $\cs=\ct_G^{\leq 0}$ it is then enough to prove that $\cs$ is closed under extensions. Therefore let $F_1,F_2\in\cs$, and let $F$ be an extension of them sitting in the distinguished triangle
\[
F_2[-1]\la F_1\la F.
\]
Assume further that, for $i=1,2$, we have distinguished triangles
\[
E_i\la F_i\la D_i,
\]
where $D_i\in\ct_G^{\leq -1}$ and $E_i\in\ct^0$. The data assembles to the diagram
\[
\xymatrix{
E_2[-1]\ar[r]& F_2[-1]\ar[r]\ar[d]& D_2[-1]\\
E_1\ar[r]& F_1\ar[r]& D_1.
}
\]
As \eqref{eq:vancomp} holds, the diagram above can be completed to a
$3\times3$ diagram, yielding the distinguished triangles
\[
E_1\la E\la E_2\qquad D_1\la D\la D_2\qquad E\la F\la D.
\]
By construction, $D\in\ct_G^{\leq -1}$ and $E\in\ct^0$. Thus (ii) in \autoref{def:wa} is automatically satisfied with $A=2$.

The argument above is essentially the same as the one in \cite[Corollary 4.3]{BV}. See also \cite[Remark 3.3]{Nee3}.
\end{remark}

Weakly approximable triangulated categories have been investigated in several papers. Let us briefly recall the main properties that will be needed in this paper:
\begin{enumerate}
\item[(A)] If $\ct$ is a weakly approximable triangulated category with weak approximation data $(G,\tau,A)$, then the \tstrs $\tau$ is in the preferred equivalence class (see \cite[Proposition 2.4]{Nee3}).

\item[(B)] Let $\ct$ be a weakly approximable triangulated category. Let $H$ be a compact generator of $\ct$ and let $\tau$ be a \tstrs in the preferred equivalence class. Then there exists an integer $B>0$ such that $(H,\tau,B)$ is a weak approximation data (see \cite[Proposition 2.6]{Nee3}).
\end{enumerate}
Both statements have parallels that work for approximable triangulated categories, see the references to \cite{Nee3} given above. Furthermore, (A) and (B) show that, once we know that a triangulated category $\ct$ is weakly approximable, then we can be very flexible in the choice of the weak approximation data. The compact generator can be changed at will, and the \tstrs can also be replaced, subject only to constraint of remaining in the preferred equivalence class. Of course the integer $A$ in the triple will have to be changed accordingly.

\begin{remark}\label{rmk:waeq}
(i) It is clear from the definition that if $\ct_1$ is a weakly approximable triangulated category, with weak approximation data $(G_1,\tau_1,A_1)$, and if $\fF\colon\ct_1\la\ct_2$ is a triangulated equivalence, then $\ct_2$ is weakly approximable as well with weak approximation data $(\fF(G_1),\fF(\tau_1),A_1)$.

(ii) In \autoref{def:wa} we stress the fact that in a weakly approximable triangulated category $\ct$ with weak approximation data $(G,\tau,A)$, the $0$-th cohomology of an object $F\in\ct^{\leq 0}$ can be approximated by $G$. The process can be iterated, as explained in \cite[Corollary 2.2]{Nee3}. More precisely, given $F\in\ct^{\leq 0}$, for any integer $m>0$, we can find a distinguished triangle
\[
E_m\la F\la D_m,
\]
where $E_m\in\ogenu G{}{1-m-A,A}$ and $D_m\in\ct^{\leq -m}$.
\end{remark}

\subsection{Examples}\label{subsec:waex}

Following the flow of the previous sections, we now want to revisit the examples from \autoref{subsec:excompgen}, in the light of the theory of weakly approximable triangulated categories.

\subsubsection*{Modules over a ring}

If $R$ is any ring, then the ring itself, seen as a compact generator for $\D(\Mod{R})$, satisfies the vanishing $\Hom(R,R[i])=0$, for all $i\geq 1$. Thus \autoref{rmk:wadep} ensures that $\D(\Mod{R})$ is weakly approximable with weak approximation data given by the triple $(R,\tau_\text{stan},1)$.

The check that $\D(\Mod{R})$ is weakly approximable can be done directly as an exercise, see \cite[3.1]{Nee3}. While carrying out the computation, one easily sees that the same triple makes $\D(\Mod{R})$ into an approximable triangulated category.

\subsubsection*{Schemes with support}

If $X$ is a quasi-compact and quasi-separated scheme, with a closed subscheme $Z\subseteq X$ such that $X\setminus Z$ is quasi-compact, then the check that $\Dqcs Z(X)$ is weakly approximable is not so easy. But it is true, as proved in

\begin{thm}[\cite{Nee2}, Theorem 3.2 (iv)]\label{thm:waschemes}
If $X$ and $Z$ are as above, then $\Dqcs Z(X)$ is weakly approximable with weak approximation data $(G,\tau_\text{stan},A)$, where $G$ is a compact generator as in \autoref{rmk:ideaproof1} and $A\gg 0$.
\end{thm}

In the special case where $Z=X$ and $X$ is separated, the category $\Dqc(X)$ is approximable (not only weakly approximable). See \cite[Lemma 3.5]{Nee3}. On the other hand, $\Dqcs Z(X)$ is not usually approximable for $X$ and $Z$ more general, but still satisfying the assumptions of \autoref{thm:waschemes}. This is explained in \cite[Remark 8.1]{Nee2}.

There is yet another geometric and interesting case which fits into the framework of \autoref{rmk:wadep}, and therefore provides another example where axiom (ii) in \autoref{def:wa} is redundant. Let $X$ be a smooth projective curve over a field. Assume further $Z=X$. It is a classical result that $\Dqc(X)$ is compactly generated by a sheaf $G$ which is a finite direct sum of ample line bundles. For a reference see, for example, \cite[Proposition 7.6]{R}. It follows that $\Hom(G,\sh[i]{G})=0$ for $i\geq 2$. Thus \autoref{rmk:wadep} immediately implies that $\Dqc(X)$ is weakly approximable.

\subsubsection*{Homotopy category of spectra}

As was observed in \cite[Example 3.2]{Nee3}, the triangulated category $\ho{\Spe}$ can be dealt with using \autoref{rmk:wadep} yet again. In \autoref{subsec:excompgen} we observed that the sphere spectrum $S^0$ is a compact generator, and $\Hom(S^0,\sh[i]{S^0})=0$ for all $i\geq 1$. Since $\tau_\text{stan}=\tau_{S^0}$,
from \autoref{rmk:wadep} we learn that the triple
$(S^0,\tau_\text{stan},2)$ provides weak approximation data for $\ho{\Spe}$. Actually one can prove more. Indeed, the triple
$(S^0,\tau_\text{stan},1)$ provides weak approximation data for $\ho{\Spe}$. Even better: it provides (strong) approximation
data for $\ho{\Spe}$.
The category $\ho{\Spe}$ is approximable, not only weakly approximable.

\subsection{The subcategories}\label{subsec:wasubcat}

Let now $\ct$ be a weakly approximable triangulated category. Fix a \tstrs $\tau=(\ct^{\leq 0},\ct^{\geq 1})$ on $\ct$ as in \autoref{def:wa}. Using this and the full subcategory of compact objects, we can define the following natural full triangulated subcategories:
\begin{itemize}
\item \emph{Bounded above objects}: $\ct^-:=\bigcup_{m=1}^\infty\ct^{\leq m}$;
\item \emph{Bounded below objects}: $\ct^+:=\bigcup_{m=1}^\infty\ct^{\geq-m}$;
\item \emph{Bounded objects}: $\ct^b:=\ct^-\cap\ct^+$,
\item \emph{Compact objects}: $\ct^c$;
\item \emph{Pseudo-compact objects}: $\ct^-_c$, where an object $F\in\ct$ belongs to $\ct^-_c$ if, for any integer $m>0$, there exists in $\ct$ a distinguished triangle $E\la F\la D$ with $E\in\ct^c$ and with $D\in\ct^{\leq-m}$. Thus
\[
\ct^-_c=\bigcap_{m=1}^\infty\big(\ct^c*\ct^{\leq-m}\big);
\]
\item \emph{Bounded pseudo-compact objects}: $\ct^b_c:=\ct^-_c\cap\ct^b$.
\item \emph{Bounded compact objects}: $\ct^{c,b}:=\ct^c\cap\ct^b=\ct^c\cap\ct^b_c$.
\end{itemize}

\begin{remark}\label{rmk:waindepsubcat}
(i) By properties (A) and (B) in \autoref{subsec:wadef}, the subcategories in the list which are defined in terms of the given \tstrs $\tau$ can be equivalently defined in terms of any \tstrs in the preferred equivalence class.

(ii) Let $\ct_1$ be a weakly approximable triangulated category and let $\fF\colon\ct_1\la\ct_2$ be a triangulated equivalence. By \autoref{rmk:waeq}(i), $\ct_2$ is weakly approximable as well. Thus both $\ct_1$ and $\ct_2$ come with the list of subcategories above. Since by (i) their definition does not depend on the choice of a \tstrs in the preferred equivalence class and, by \autoref{rmk:preststr2}, $\fF$ preserves the preferred equivalence class, the equivalence $\fF$ restricts to equivalences between the corresponding subcategories of $\ct_1$ and $\ct_2$.
\end{remark}

By their explicit definitions, the subcategories listed above sit in the following commutative diagram:
\begin{equation}\label{eq:incl}
\xymatrix@C+40pt@R-5pt{
&\ct&\\
\ct^-\ar@{^{(}->}[ur]&\ct^b\ar@{^{(}->}[u]\ar@{_{(}->}[l]\ar@{^{(}->}[r]&\ct^+\ar@{_{(}->}[ul]\\
\ct^-_c\ar@{^{(}->}[u]&\ct^b_c\ar@{^{(}->}[u]\ar@{_{(}->}[l]&\\
\ct^c\ar@{^{(}->}[u]&\ct^{c,b}.\ar@{^{(}.>}[u]\ar@{_{(}->}[l]
 \ar@/_2pc/@{^{(}->}[uu]&
}
\end{equation}
It is worth noting that the above-listed categories turn out to be  natural, classical old friends in two out of the three classes of examples described in \autoref{subsec:waex}:

\medskip

\begin{center}
	\begin{tabular}{|c|c|c|}

		\hline
		& $Z\subseteq X$ as in \autoref{subsec:excompgen}  & $R$ a ring\\
		\hline
		$\ct$ & $\Dqcs Z(X)$ & $\D(\Mod{R})$\\
		$\ct^-$ & $\Dqcmis Z(X)$ & $\D^-(\Mod{R})$\\
		$\ct^+$ & $\Dqcpls Z(X)$ & $\D^+(\Mod{R})$\\
		$\ct^b$ & $\Dqcbs Z(X)$ & $\D^b(\Mod{R})$\\
		$\ct^c$ & $\dperfs ZX$ & $\dperf{R}$\\
		$\ct^-_c$ & $\Dqcps Z(X)$ & $\K^-(\proj{R})$\\
		$\ct^b_c$ & $\Dqcpbs Z(X)$ & $\K^{-,b}(\proj{R})$\\
  $\ct^{c,b}$ & $\dperfs ZX$ & $\dperf{R}$\\
		\hline
	\end{tabular}
\end{center}

\medskip

Comparing the table with the diagram \eqref{eq:incl} it is clear that, in these two examples, the subcategories $\ct^{c,b}$ and $\ct^c$ coincide. But this is definitely not always true. In the special case $\ct=\ho{\Spe}$, we have $\ct^{c,b}=\ct^b_c\cap\ct^c=\{0\}$ whereas $\ct^{b}_c\ne\{0\}\ne\ct^c$.

The main result of this section proves that all the solid inclusions $\ca\mono\cb$ in \eqref{eq:incl} are intrinsic. Roughly speaking, this means that there is a recipe, depending only on the triangulated structure on $\cb$, that describes which objects in $\cb$ belong to the full subcategory $\ca$. In order to deal also with the inclusion represented by a dotted arrow, we need to restrict to a special class of weakly approximable triangulated categories.

\begin{definition}\label{def:coherent}
A weakly approximable triangulated category $\ct$ is \emph{coherent} if, for any \tstrs $\tst\ct$ in the preferred equivalence class, there
exists an integer $N>0$ such that every
object $Y\in\ct^-_c$ admits a distinguished triangle
$X\la Y\la Z$ with
$X\in\ct^-_c\cap\ct^{\leq N}$ and with $Z\in\ct^b_c\cap\ct^{\geq0}$.
\end{definition}



We are now ready to state the main result of this section.

\begin{thm}[\cite{CNS2}, Theorem B]\label{thm:main1}
All the inclusions in diagram \eqref{eq:incl}, represented by the solid arrows $\ca\mono\cb$, are \emph{invariant under triangulated equivalences}. By this we mean that given a pair of weakly approximable triangulated
categories $\ct,\ct'$, as well as matching inclusions
$\ca\mono\cb\mono\ct$ and
$\ca'\mono\cb'\mono\ct'$ from diagram \eqref{eq:incl},
then any triangulated equivalence $\cb\la\cb'$ must restrict to a triangulated equivalence $\ca\la\ca'$.

Furthermore, the same is true for the unique inclusion in diagram \eqref{eq:incl} represented by a dotted arrow, provided we further assume that one of the two conditions below holds.
\begin{itemize}
\item $\ct,\ct'$ are coherent, or
\item $\ct^c\subseteq\ct^b_c$, ${\ct'}^c\subseteq{\ct'}^b_c$ and
$^\perp(\ct^b_c)\cap\ct^-_c=\{0\}={^\perp({\ct'}^b_c)}\cap{\ct'}^-_c$.
\end{itemize}
\end{thm}

It is clear that, if $\cb=\ct$ and $\cb'=\ct'$, then the result is just \autoref{rmk:waindepsubcat}(ii). The hard part of the proof consists in dealing with the other cases. Furthermore, at the moment, it is not clear whether the assumptions in the second part of the statement are really needed in order to deal with the dotted inclusion. It is worth recalling that, given a full subcategory $\cs\subseteq\ct$, the collection of objects of the full subcategory $^\perp\cs\subseteq\ct$ is $\{A\in\ct:\Hom(A,S)=0\text{ for all }S\in\cs\}$.

As for the applications, it is useful to note the following properties.
\begin{itemize}
\item If $R$ is a left coherent ring, then it is easy to prove that $\D(\Mod{R})$ is coherent as well. See \cite[Example 10.1.2]{CNS2} for a complete discussion.
\item If $\ct$ is coherent and $\ct^c\subseteq\ct^b_c$, then $^\perp(\ct^b_c)\cap\ct^-_c=\{0\}$. See \cite[Proposition 10.1.3]{CNS2}.
\item The condition $\ct^c\subseteq\ct^b_c$ and $^\perp(\ct^b_c)\cap\ct^-_c=\{0\}$ holds if $\ct$ is either $\D(\Mod{R})$ with $R$ a commutative ring, or $\Dqcs Z(X)$ with $X$ and $Z$ as before. See \cite[Proposition 10.2.1]{CNS2}.
\end{itemize}



\section{Enhancements}\label{sec:enh}

In this section we change perspective. The possible higher categorical incarnations, of the triangulated categories discussed in the previous sections, come to play a key role. In \autoref{subsec:enhdefex} and \autoref{subsec:dgex} we recall the basic definitions, properties and examples, while in \autoref{subsec:enhuniq} we move on to new results. We will survey developments concerning the uniqueness of enhancements. This will be used, in the last section, to discuss the so called Morita theory for schemes. In \autoref{subsec:enhwip} we mention some ongoing work, still in progress.

\subsection{Definitions and first constructions}\label{subsec:enhdefex}

We should begin this section with an admission. But first the background: one of the crucial observation in the field is that most triangulated categories, of algebro-geometric interest, admit higher categorical enhancements. And in many important cases the enhancement is unique, in a sense that can be made precise.

That said, enhancements come in more than one flavor, and there is no consensus yet on which of the variants is best. As it turns out, our newest results are largely flavor-independent, they hold for the entire spectrum of different enhancement types. And now we come to our admission: because of this robustness of the
results, of the fact that the latest theorems apply across the spectrum of
enhancement types, on occasion, in this survey we will allow ourselves to be imprecise about which explicit enhancement type we adopt. 

There are at least three enhancement types which have been used extensively in the recent literature: differential grade (dg) categories, $A_\infty$--categories and $\infty$-categories. In this section, we restrict attention to the first of these. We work with dg categories, largely because our recent (published) paper, \cite{CNS1}, which this section surveys, uses the dg machinery. There are some comments in \cite{CNS1} about the extent to which the results remain valid for other enhancement types, and there is literature studying the issue. The interested reader is referred to \cite{COS1,COS2} for the comparison between dg and $A_\infty$--enhancements, and to \cite[Theorem 14]{A1}, as well as \cite{Co} and \cite{Do}, for a comparison of $\infty$-categorical enhancements with dg and $A_\infty$--enhancements.

Let $\kk$ be a commutative ring. We begin with the following basic definitions.

\begin{definition}\label{def:dg}
(i) A \emph{dg category} is a $\kk$-linear category $\cc$ whose morphism spaces $\Hom_\cc\left(A,B\right)$ are complexes of $\kk$-modules and the composition maps $\Hom_\cc(B,C)\otimes_{\kk}\Hom_\cc(A,B)\to\Hom_\cc(A,C)$ are morphisms of complexes, for all $A,B,C$ in $\cc$.

(ii) A \emph{dg functor} $\fF\colon\cc_1\to\cc_2$ between two dg categories is a $\kk$-linear functor such that the maps $\Hom_{\cc_1}(A,B)\to\Hom_{\cc_2}(\fF(A),\fF(B))$ are morphisms of complexes, for all $A,B$ in $\cc_1$.
\end{definition}

Given a dg category $\cc$, we can consider its \emph{homotopy category} $H^0(\cc)$ with the same objects of $\cc$ and such that $\Hom_{H^0(\cc)}(A,B):=H^0(\Hom_\cc(A,B))$), for all $A,B$ in $\cc$. Clearly a dg functor $\fF\colon\cc_1\to\cc_2$ induces a $\kk$-linear functor $H^0(\fF)\colon H^0(\cc_1)\la H^0(\cc_2)$ between the corresponding homotopy categories. A dg functor $\fF$ is a \emph{quasi-equivalence} if the maps $\Hom_{\cc_1}(A,B)\to\Hom_{\cc_2}(\fF(A),\fF(B))$ are quasi-isomorphisms, for all $A,B$ in $\cc_1$, and $H^0(\fF)$ is an equivalence.

Let $\dgCat$ be the category of all (small) dg categories. One can then localize $\dgCat$ by inverting all quasi-equivalences. We denote by $\Hqe$ such a localization. In this section, this is the category where all models of triangulated categories will be compared.

\begin{remark}\label{rmk:dgcompAinfty}
If $\cc_1$ and $\cc_2$ are dg categories and $f\colon\cc_1\la\cc_2$ is a morphism in $\Hqe$, then we get an induced functor $H^0(f)\colon H^0(\cc_1)\la H^0(\cc_2)$, which is well defined up to isomorphism. If $f$ is an isomorphism in $\Hqe$, then $H^0(f)$ is an equivalence.
\end{remark}

As dg functors between two dg categories $\cc_1$ and $\cc_2$ form a dg category $\dgHom(\cc_1,\cc_2)$, for every dg category $\cc$ we can construct a much larger one
\[
\dgMod{\cc}:=\dgHom(\cc\opp,\dgC(\Mod{\kk}))
\]
whose objects are called \emph{(right) dg $\cc$-modules} (see \autoref{subsec:dgex} for the definition of the dg category $\dgC(\Mod{\kk})$). The category $H^0(\dgMod{\cc})$ has a natural triangulated structure (see \cite[Section 2.2]{K3}). On the other hand, given a dg category $\cc$, the map defined on objects by $A\mapsto\Hom_\cc(\farg,A)$ extends to a fully faithful dg functor $\dgYon\colon\cc\to\dgMod{\cc}$ called the \emph{dg Yoneda embedding}.

\begin{definition}\label{def:dgpreptr}
(i) A dg category $\cc$ is \emph{pretriangulated} if the essential image of $H^0(\dgYon)$ is a full triangulated subcategory of $H^0(\dgMod{\cc})$.

(ii) Given a triangulated category $\ct$ an \emph{enhancement} of $\ct$ is a pair $(\cc,\fF)$, where $\cc$ is a pretriangulated dg category and $\fF\colon H^0(\cc)\la\ct$ is a triangulated equivalence.

(iii) A triangulated category $\ct$ is \emph{algebraic} if it has an enhancement.
\end{definition}

One can additionally show that the image of $\dgYon$ is contained in the full dg subcategory $\hproj{\cc}$ of $\dgMod{\cc}$ whose objects are \emph{h-projective} dg $\cc$-modules. More precisely, an object $M\in\dgMod{\cc}$ is \emph{h-projective} if $\Hom_{H^0(\dgMod{\cc})}(M,N)=0$ for every $N\in\dgAc{\cc}$. Here $\dgAc{\cc}$ is the full dg subcategory of $\dgMod{\cc}$ whose objects are acyclic, meaning that $N(A)$ is an acyclic complex for every $A\in\cc$.

\begin{remark}\label{rmk:dgimp}
(i) Let $\cc$ be a pretriangulated dg category such that the triangulated category $H^0(\cc)$ has coproducts, and assume that $\cb\subseteq\cc$ is a full dg subcategory whose objects form a set of compact generators for $H^0(\cc)$. Then, by \cite[Proposition 1.17]{LO}, there is a natural isomorphism $\hproj{\cb}\cong\cc$ in $\Hqe$.

(ii) Every (iso)morphism $f\colon\cc_1\la\cc_2$ in $\Hqe$ extends to a unique (iso)morphism $\tilde f\colon\hproj{\cc_1}\la\hproj{\cc_2}$ in $\Hqe$, which is compatible with the Yoneda embeddings and such that $H^0(\tilde f)$ preserves coproducts (see \cite{Dr} for more details).
\end{remark}

\subsection{Examples}\label{subsec:dgex}

Let us now discuss few basic examples which are relevant for our discussion.

Assume that $\ca$ is a $\kk$-linear additive category. We can form the dg category $\dgC(\ca)$, whose objects are (unbounded) complexes of objects in $\ca$. As graded modules, morphisms are defined as
\[
\Hom_{\dgC(\ca)}(A,B)^n:=\prod_{i\in\ZZ}\Hom_\ca(A^i,B^{n+i})
\]
for every $A,B\in\dgC(\ca)$ and for every $n\in\ZZ$. The composition of morphisms is the obvious one while the differential is defined on a homogeneous element $f\in\Hom_{\dgC(\ca)}(A,B)^n$ by $d(f):=d_B\comp f-(-1)^nf\comp d_A$.

It is easy to see that $\dgC(\ca)$ is pretriangulated and there is a natural triangulated equivalence $H^0(\dgC(\ca))\cong\K(\ca)$. Thus the triangulated category $\K(\ca)$ is algebraic.

Assume now that $\ca$ is abelian. As any full triangulated subcategory of an algebraic triangulated category is algebraic as well, the full subcategory $\K_\text{acy}(\ca)\subseteq\K(\ca)$ consisting of acyclic complexes is algebraic. Indeed, it is realized as the homotopy category of the full dg subcategory $\dgCacy(\ca)\subseteq\dgC(\ca)$ consisting of all acyclic complexes.

In order to provide an enhancement for $\D(\ca)$, which is the Verdier quotient $\K(\ca)/\K_\text{acy}(\ca)$ by definition, we would like to be able to perform a similar quotient at the dg level. Such a construction goes under the name of \emph{Drinfeld quotient} and can be summarized as follows.

Given a dg category $\cc$, one can produce a h-flat dg category $\cc^\text{hf}$ with a quasi-equivalence $\cc^\text{hf}\la\cc$ (hence $\cc$ and $\cc^\text{hf}$ are isomorphic in $\Hqe$). Recall that a dg category is h-flat if its complexes of morphisms are homotopically flat over $\kk$, meaning that tensoring with them preserves acyclic complexes. Note that the construction of $\cc^\text{hf}$ in \cite[Appendix B]{Dr} can be modified and made functorial by \cite[Proposition 3.10]{CNS1}. Furthermore, if $\cb\subseteq\cc$ is a full dg subcategory, we get the corresponding h-flat full dg subcategory $\cb^\text{hf}\subseteq\cc^\text{hf}$. Then we call the \emph{Drinfeld quotient} of $\cc$ by $\cb$ and denote by $\cc/\cb$ what should be denoted more appropriately by $\cc^\text{hf}/\cb^\text{hf}$. We refer to \cite[Section 3.1]{Dr} for the details of the construction. The remarkable property is that, if $\cc$ and $\cb$ are pretriangulated, then $\cc/\cb$ is pretriangulated as well and there is a natural triangulated equivalence
\begin{equation}\label{eq:dgcoh}
H^0\left(\cc/\cb\right)\cong H^0(\cc)/H^0(\cb).
\end{equation}
In the special case of an abelian category $\ca$, we consider the dg category
\[
\Ddg(\ca):=\dgC(\ca)/\dgCacy(\ca).
\]
By \eqref{eq:dgcoh}, we have a natural triangulated equivalence
\[
\D(\ca)\cong H^0(\Ddg(\ca))
\]
and thus $\D(\ca)$ is algebraic. As full triangulated subcategories of algebraic ones are algebraic, we immediately get the following:

\begin{prop}\label{prop:dgalg}
The following triangulated categories are algebraic:
\begin{itemize}
\item $\Ka(\ca)$ when $\ca$ is an additive category and $\D^?(\ca)$ when $\ca$ is an abelian category, for $?=b,+,-,\emptyset$;
\item $\D(\Mod{R})$ and all its subcategories from \autoref{subsec:wasubcat}, for any ring $R$;
\item $\Dqcs Z(X)$ and all its subcategories from \autoref{subsec:wasubcat}, for $X$ a quasi-compact and quasi-separated scheme with a closed subscheme $Z$ such that $X\setminus Z$ is quasi-compact.
\end{itemize}
\end{prop}

\begin{remark}\label{rmk:dgnonalg}
It is well-known that $\ho{\Spe}$ is not algebraic---see, for example, \cite[Proposition 3.6]{CS:surv1}. On the other hand it has natural $\infty$-categorical models.
\end{remark}

\subsection{Uniqueness}\label{subsec:enhuniq}

Coming back to the triangulated categories in \autoref{prop:dgalg}, which we know to be algebraic, we want to investigate the uniqueness of their enhancements. To make this rigorous, we start with the following definition:

\begin{definition}\label{def:dguniq}
Let $\ct$ be an algebraic triangulated category.

(i) The category $\ct$ \emph{has a unique enhancement} if, given two enhancements $(\cc_1,\fF_1)$ and $(\cc_2,\fF_2)$ of $\ct$, the dg categories $\cc_1$ and $\cc_2$ are isomorphic in $\Hqe$.

(ii) The category $\ct$ \emph{has a strongly unique enhancement} if, given two enhancements $(\cc_1,\fF_1)$ and $(\cc_2,\fF_2)$ of $\ct$, there exists an isomorphism  $f\colon\cc_1\la \cc_2$ in $\Hqe$ such that $\fF_2\circ H^0(f)=\fF_1$.
\end{definition}

The study of the uniqueness of enhancements for triangulated categories has a long history. It was motivated by the desire to show that constructions, which appeal to enhancements, are independent of the enhancements chosen. This is also the case in the current paper, in a sense to be elaborated below.

We begin with
a trivial observation.

\begin{remark}\label{rmk:dgimplift}
Assume that $\ct_1$ and $\ct_2$ are algebraic triangulated categories, and let $\fF\colon\ct_1\la\ct_2$ be a $k$-linear triangulated equivalence. If $\ct_1$ has a unique enhancement then so does $\ct_2$, and any enhancement $\mathbf{T}_1$ of
$\ct_1$ is isomorphic in $\Hqe$ to any enhancement $\mathbf{T}_2$ of $\ct_2$.
Of course, applying $H^0$ to the isomorphism
$\mathbf{T}_1\cong\mathbf{T}_2$ gives an equivalence
$\ct_1\cong\ct_2$. This equivalence need
not agree with the given $\fF$.

For us the important advantage, of having a triangulated equivalence coming from an isomorphism in $\Hqe$, is that it extends. We have in a sense already seen one example: by \autoref{rmk:dgimp}, if $\cc_1$ and $\cc_2$ are pretriangulated dg categories and $\cb_i\subseteq\cc_i$ is the full dg subcategory such that $H^0(\cb_i)=H^0(\cc_i)^c$, for $i=1,2$, then an isomorphism $\cb_1\cong\cb_2$ in $\Hqe$ extends in $\Hqe$ to an isomorphism $\cc_1\cong\cc_2$. 
\end{remark}

With this in mind, we are ready to state the second main result of this paper.

\begin{thm}[\cite{CNS2}, Theorems A and B]\label{thm:main2}
{\rm (i)} If $\ca$ is an abelian category, then $\D^?(\ca)$ has a unique enhancement, for $?=b,+,-,\emptyset$.

{\rm (ii)} Let $X$ be a quasi-compact and quasi-separated scheme. Then the categories $\Dqca(X)$ and $\dperf{X}$ have a unique dg enhancement, for $?=b,+,-,\emptyset$.
\end{thm}

\begin{ex}\label{ex:dgmod}
The first part of the result above, applied to the special case $\ca=\Mod{R}$, proves that $\D(\Mod{R})$ has a unique enhancement. But this can easily be proved directly, and was done much earlier. It boils down to a simple computation with generators. The same applies if one wants to prove uniqueness of enhancements for $\K^b(\proj R)=\D(\Mod{R})^c$. See \cite[Proposition 2.6]{LO} for more details.
\end{ex}

\autoref{thm:main2} is, for now, the most general result in the literature. It does have several predecessors: one of the first papers, where the problem of uniqueness appeared, was \cite{BLL}---where it was conjectured that, for a smooth projective scheme, the bounded derived category of coherent sheaves should have a unique enhancement. The conjecture was proved in the groundbreaking paper \cite{LO}, and the results therein were extended in \cite{CS6}. Part of the results in the latter paper have been recently reproved in \cite{GR}. Finally the paper \cite{A1} brought to the theory new results and conjectures by resetting the problem in the context of stable $\infty$-categories.

As for strong uniqueness of enhancements, the first results were obtained in \cite{LO}. Despite some progress in the following years (see, for instance, \cite{CS1}, \cite{Ol} and \cite{Lo}), so far proving that a given triangulated category has a strongly unique enhancement remains a much more challenging problem than proving uniqueness. In particular, if we restrict to the triangulated categories in the statement of \autoref{thm:main2}, we can say that strong uniqueness is known to hold only in some special cases, and is still a difficult open problem in general.

\subsection{Work in progress}\label{subsec:enhwip}

In an ongoing research project, we plan to combine the theory of weakly approximable triangulated categories with techniques coming from higher categories. There are actually two directions which are very promising and which might lead to significant improvements of \autoref{thm:main2} and which might be applicable to the discussion in the next section.

We first observe that, by \autoref{rmk:dgnonalg}, we have examples of weakly approximable triangulated categories which are not algebraic. On the other hand, if $\ct$ is weakly approximable, we can look at all its natural subcategories in \eqref{eq:incl} and even if $\ct$ is not algebraic, it might (in theory)\footnote{We don't know any examples of this.} happen that one of its subcategories $\cd\mono\ct$ in \eqref{eq:incl} is algebraic. In this case, we already explained that all the subcategories of $\cd$ are automatically algebraic.

From this perspective, a natural problem can be formulated as follows:

\begin{pb}\label{pb:progr1}
  Let $\ct$ be a weakly approximable triangulated category and let $\ca\mono\cb$ be one of the natural inclusions in \eqref{eq:incl}. Assume that $\cb$ is algebraic. Determine when and to which extent the (strong) uniqueness of enhancement of $\ca$ implies the (strong) uniqueness of enhancement of $\cb$, and
  vice versa.
\end{pb}

A result in this direction would not be surprising. If $\cb=\ct$ and $\ca=\ct^c$, then $\ca$ having a unique enhancement implies that so does $\cb$. This is a direct consequence of \autoref{rmk:dgimp}. But this is, so far, the only known result in this direction. To prove that uniqueness of enhancement for
$\ca$ implies uniqueness of enhancement for $\cb$ would, presumably, require a reconstruction theorem, in the style of \autoref{rmk:dgimp}(i), for all subcategories. This is clearly more challenging and we will address \autoref{pb:progr1} in a forthcoming paper.

Hopefully, once this is achieved, we would get a wealth of new results already in the geometric and algebraic cases. The idea would be to combine, whatever results we can prove towards \autoref{pb:progr1}, with what we already know from  \autoref{thm:main2}. In particular, let $R$ be a ring and let $X$ be a quasi-compact quasi-separated scheme. If we are going to be optimistic and assume the answer to \autoref{pb:progr1} is a resounding Yes in all cases, then  the following triangulated categories would have a unique enhancement:
\begin{itemize}
\item $\K^-(\proj R)$ and $\K^{-,b}(\proj R)$;
\item $\Dqcp(X)$ and $\D_\scat{qc}^{p,b}(X)$.
\end{itemize}

The careful reader might have noticed that all the geometric categories, for which we can prove uniqueness of enhancement, do not involve support. Whereas the results on weakly approximable triangulated categories, which we have seen in the previous sections, all hold in the relative case. The second part of our ongoing work in progress aims to find a remedy for this, by addressing the following.

\begin{pb}\label{pb:progr2}
Let $X$ be a quasi-compact and quasi-separated scheme with a closed subscheme $Z$ such that $X\setminus Z$ is quasi-compact. Show that $\dperfs ZX$ has a unique enhancement. 
\end{pb}

We will discuss in the next section an interesting application of a result in this direction. In any case it should be clear that, a positive solution to both \autoref{pb:progr1} and \autoref{pb:progr2}, would immediately lead to a sharp improvement of the second part of \autoref{thm:main2}. Given a quasi-compact and quasi-separated scheme $X$ with a closed subscheme $Z$ such that $X\setminus Z$ is quasi-compact, it would follow that the triangulated categories
\[
\Dqcsa Z(X)\qquad\Dqcps Z(X)\qquad
\Dqcpbs {Z}(X),
\]
for  $?=b,+,-,\emptyset$, would also have unique enhancements.

\section{Morita theory and weakly approximable triangulated categories}\label{sec:Morita}

We conclude this survey by explaining a result, which extends and vastly generalizes all known theorems about derived invariance for rings. As we will see in \autoref{subsec:Rickard} such a result can be thought of as a Morita-type theorem for schemes. What makes it the natural conclusion of this survey is that its proof is already complete, and is a place where our results on weakly approximable categories and their subcategories, and on uniqueness of enhancements, fruitfully meet. In \autoref{subsec:appl} we provide some additional geometric applications.

\subsection{Rickard's theorem for schemes}\label{subsec:Rickard}

In this section we outline one major application of our previous results \autoref{thm:main1} and \autoref{thm:main2}, and of its potential extensions discussed in the previous section.

The starting point of our investigation is the following beautiful result by Rickard.

\begin{thm}[\cite{Rickard91,Rickard89b}]\label{thm:Rickard}
Let $R$ and $S$ be two rings. Then
the following are
equivalent:
\begin{enumerate}
\item
There exists a triangulated equivalence $\D(\Mod R)\cong\D(\Mod S)$.
\item
There exists a triangulated equivalence $\D^-(\Mod R)\cong\D^-(\Mod S)$.
\item
There exists a triangulated equivalence $\D^+(\Mod R)\cong\D^+(\Mod S)$.
\item
There exists a triangulated equivalence $\D^b(\Mod R)\cong\D^b(\Mod S)$.
\item
There exists a triangulated equivalence $\K^-(\proj R)\cong\K^-(\proj S)$.
\item
There exists a triangulated equivalence $\K^b(\proj R)\cong\K^b(\proj S)$.
\setcounter{enumiv}{\value{enumi}}
\end{enumerate}
Moreover the six equivalent conditions above imply:
\begin{enumerate}
\setcounter{enumi}{\value{enumiv}}
\item
There exists a triangulated equivalence $\K^{-,b}(\proj R)\cong\K^{-,b}(\proj S)$.
\end{enumerate}
The converse holds if we assume further that $R$ and $S$ are both left coherent, in which case the equivalence in (7) can be rewritten as $\D^b(\mmod R)\cong\D^b(\mmod S)$.
\end{thm}

The beauty of this result is that, as a consequence, any of the triangulated categories associated to a ring, in (1)--(6), carries the exact same information as any of the others. And this can be extended to (7) in the left coherent case. It makes sense to formulate the following definition.

\begin{definition}\label{def:dereqrings}
Two rings $R$ and $S$ are \emph{derived equivalent} if any of the conditions (1)--(6) is satisfied.
\end{definition}

In the language of \cite{Rickard91,Rickard89b}, \autoref{thm:Rickard} is an instance of the so-called \emph{Morita theory} for derived categories of rings. The reason is that, by analogy with classical Morita theory, Rickard proves that the equivalences in the list are produced by tensoring with a suitable complex of bimodules. It is interesting to note that this implies, in particular, that derived equivalence is left/right symmetric, in the sense that two rings $R$ and $S$ are derived equivalent if and only if $R\opp$ and $S\opp$ are.

The discussion in the previous sections might have shown the reader that there should be nothing special about $\D(\Mod R)$, and that we should expect a result similar to \autoref{thm:Rickard} in a much more general setting. In fact, \autoref{thm:main1} in combination with the argument of \autoref{rmk:dgimp} easily yields the following major generalization of Rickard's result.

\begin{thm}\label{thm:main3}
  Let $\ct$ and $\ct'$ be two weakly approximable triangulated categories. Assume that
\begin{itemize}
\item
$\ct$ and $\ct'$ are algebraic.
\item  
$\ct^c$ or ${\ct'}^c$ has a unique enhancement.
\end{itemize}
Then the following are equivalent:
\begin{enumerate}
\item There exists a triangulated equivalence $\ct\cong\ct'$.
\item There exists a triangulated equivalence $\ct^+\cong{\ct'}^+$.
\item There exists a triangulated equivalence $\ct^-\cong{\ct'}^-$.
\item There exists a triangulated equivalence $\ct^b\cong{\ct'}^b$.
\item There exists a triangulated equivalence $\ct^-_c\cong{\ct'}^-_c$.
\item There exists a triangulated equivalence $\ct^c\cong{\ct'}^c$.
\setcounter{enumiv}{\value{enumi}}
\end{enumerate}
Moreover the six equivalent conditions above imply:
\begin{enumerate}
\setcounter{enumi}{\value{enumiv}}
\item There exists a triangulated equivalence $\ct^b_c\cong{\ct'}^b_c$.
\end{enumerate}
The converse holds if we assume further that $\ct^c\subseteq\ct^b_c$, ${\ct'}^c\subseteq{\ct'}^b_c$ and $^\perp(\ct^b_c)\cap\ct^-_c=\{0\}={^\perp({\ct'}^b_c)}\cap{\ct'}^-_c$.
\end{thm}

By the discussion after \autoref{thm:main1} it should be clear that, in the special case where $\ct=\D(\Mod R)$ and $\ct'=\D(\Mod S)$, \autoref{thm:main3} delivers \autoref{thm:Rickard}. Actually a slightly stronger result, since in the last part $R$ and $S$ can be assumed to be either left coherent or commutative. Note that $\K^b(\proj R)=\D(\Mod{R})^c$ has a unique enhancement by \autoref{ex:dgmod}. More importantly, using also \autoref{thm:main2}, in the geometric case we obtain the following surprising, new result.

\begin{cor}\label{thm:Morita}
Let $X$ and $X'$ be quasi-compact and quasi-separated schemes with a closed subscheme $Z'\subseteq X'$ such that $X'\setminus Z'$ is quasi-compact. Then the following are equivalent:
\begin{enumerate}
\item There exists a triangulated equivalence $\Dqc(X)\cong\Dqcs{Z'}(X')$.
\item There exists a triangulated equivalence $\Dqcpl(X)\cong\Dqcpls{Z'}(X')$.
\item There exists a triangulated equivalence $\Dqcmi(X)\cong\Dqcmis{Z'}(X')$.
\item There exists a triangulated equivalence $\Dqcb(X)\cong\Dqcbs{Z'}(X')$.
\item There exists a triangulated equivalence $\Dqcp(X)\cong\Dqcps{Z'}(X')$.
\item There exists a triangulated equivalence $\dperf X\cong\dperfs{Z'}{X'}$.
\item There exists a triangulated equivalence $\Dqcpb(X)\cong\Dqcpbs{Z'}(X')$.
\end{enumerate}
\end{cor}

Thus it follows that all the a priori different triangulated categories of the left-hand-sides in (1)--(7) carry again the exact same information. In analogy with \autoref{def:dereqrings} we can formulate the following:

\begin{definition}\label{def:dereqschemes}
Two quasi-compact and quasi-separated schemes $X$ and $X'$ are \emph{derived equivalent} if any of the conditions (1)--(7) above (with $Z'=X'$) is satisfied.
\end{definition}

\begin{remark}\label{rmk:proof}
(i) If two quasi-compact and quasi-separated schemes $X$ and $X'$ are derived equivalent then, thanks to uniqueness of enhancement for $\Dqc(X)$, one can always find a triangulated equivalence $\Dqc(X)\cong\Dqc(X')$ which is induced by an isomorphism in $\Hqe$. By the seminal work of To\"en \cite{To} (see also \cite{COS1,COS2,CS} for different proofs) such an isomorphism is induced by a suitable dg bimodule. Thus \autoref{thm:Morita} can be thought of as a Morita-type theorem for schemes.

(ii) Assuming \autoref{pb:progr2} turns out to have a positive answer, we can improve \autoref{thm:Morita} by replacing all the triangulated categories involving the scheme $X$ with relative versions for a given closed subscheme $Z\subseteq X$ such that $X\setminus Z$ is quasi-compact.
\end{remark}

\subsection{Applications}\label{subsec:appl}

If $X$ is a quasi-compact and quasi-separated scheme, the natural inclusion  $\dperf{X}\mono \D^{p,b}(X)$ is not always an equality. Thus the following makes sense.

\begin{definition}\label{def:catsing}
If $X$ is a quasi-compact and quasi-separated scheme, then the Verdier quotient
\[
\D_\scat{sing}(X):=\D^{p,b}(X)/\dperf{X}
\]
is the \emph{category of singularities} of $X$.
\end{definition}

Assuming that $X$ is noetherian we have an identification $\D^{p,b}(X)=\D^b(\scat{Coh}(X))$, the above simplifies to the standard category of singularities. Still assuming that
$X$ is noetherian, one defines the \emph{stable derived category} of $X$ to be
\[
\scat{S}(\scat{Qcoh}(X)):=\K_\text{acy}(\mathrm{Inj}(\scat{Qcoh}(X))\ ,
\]
that is the homotopy category of acyclic complexes of injectives in $\scat{Qcoh}(X)$. It is known that $\scat{S}(\scat{Qcoh}(X))$ is compactly generated, and $\scat{S}(\scat{Qcoh}(X))^c=\D_\scat{sing}(X)$. As a straightforward application of \autoref{thm:Morita} and \autoref{rmk:dgimp}, we obtain the following result.

\begin{cor}\label{cor:singstab}
Let $X$ and $X'$ be a quasi-compact and quasi-separated schemes which are derived equivalent. Then
\begin{itemize}
\item[{\rm (i)}] There is a triangulated equivalence $\D_\scat{sing}(X)\cong \D_\scat{sing}(X')$;
\item[{\rm (ii)}] If $X$ and $X'$ are also noetherian, then there is a triangulated equivalence $\scat{S}(\scat{Qcoh}(X))\cong \scat{S}(\scat{Qcoh}(X'))$.
\end{itemize}
\end{cor}

Recall that a noetherian scheme $X$ is regular if and only if the inclusion $\dperf{X}\mono\D^{p,b}(X)$ is an equality, that is if and only if $\D_\scat{sing}(X)$ is trivial.  We deduce the following result.

\begin{cor}\label{cor:re}
Let $X$ and $X'$ be noetherian schemes which are derived equivalent. Then $X$ is regular if and only if $X'$ is regular.
\end{cor}

\begin{remark}\label{rmk:noderinv}
The above result shows that being regular is a derived invariant for noetherian schemes. Of course one can ask if other interesting properties are derived invariant. For instance, it seems to be unknown if being noetherian is a derived invariant for quasi-compact and quasi-separated schemes. Note however that there is an  unpublished example by Rickard proving that being left noetherian (or left coherent) is not a derived invariant for (noncommutative) rings.
\end{remark}


\bigskip

{\small\noindent{\bf Acknowledgements.} Part of this work was carried out while the third author was visiting the Institut des Hautes \'Etudes Scientifiques (Paris) whose warm hospitality is gratefully acknowledged.}



\begin{thebibliography}{99}

\bibitem{Alonso-Jeremias-Souto03} L.\ Alonso~Tarr{\'{\i}}o, A.\ Jerem{\'{\i}}as~L{\'o}pez, M.~J.\ Souto~Salorio,\emph{Construction of {$t$}-structures and equivalences of derived categories}, Trans.\ Amer.\ Math.\ Soc.\ {\bf 355} (2003), 2523--2543.

\bibitem{A1} B.\ Antieau, \emph{On the uniqueness of infinity-categorical enhancements of triangulated categories}, arXiv:1812.01526.

\bibitem{BLL} A.\ Bondal, M.\ Larsen, V.\ Lunts, \emph{Grothendieck ring of pretriangulated categories}, Int.\ Math.\ Res.\ Not.\ {\bf 29} (2004), 1461--1495.

\bibitem{BV} M.V.\ Bondarko, S.V.\ Vostokov, \emph{On weakly negative subcategories, weight structures, and (weakly) approximable triangulated categories}, Lobachevskii J.\ Math.\ {\bf 41} (2020), 151--159.

\bibitem{CNS1} A.\ Canonaco, A.\ Neeman, P.\ Stellari, \emph{Uniqueness of enhancements for derived and geometric categories}, Forum Math.\ Sigma {\bf 10} (2022), 1--65.

\bibitem{CNS2} A.\ Canonaco, A.\ Neeman, P.\ Stellari, \emph{The passage among the subcategories of weakly approximable triangulated categories}, (Appendix by C.\ Haesemeyer), arXiv:2402.04605.
	
\bibitem{COS1} A.\ Canonaco, M.\ Ornaghi, P.\ Stellari, \emph{Localizations of the categories of $A_\infty$-categories and Internal Homs}, Doc. Math. {\bf 24}, (2019) 2463-2492.

\bibitem{COS2} A.\ Canonaco, M.\ Ornaghi, P.\ Stellari, \emph{Localizations of the categories of $A_\infty$-categories and Internal Homs over a ring}, arXiv:2404.06610.

\bibitem{CS1} A.\ Canonaco, P.\ Stellari, \emph{Fourier-Mukai functors in the supported case}, Compositio Math.\ {\bf 150} (2014), 1349--1383.

\bibitem{CS} A.\ Canonaco, P.\ Stellari, \emph{Internal Homs via extensions of dg functors}, Adv.\ Math.\ {\bf 277} (2015), 100--123.

\bibitem{CS:surv1} A.\ Canonaco, P.\ Stellari, \emph{A tour about existence and uniqueness of dg enhancements and lifts}, J.\ Geom.\ Phys.\ {\bf 122} (2017), 28--52.

\bibitem{CS6} A.\ Canonaco, P.\ Stellari, \emph{Uniqueness of dg enhancements for the derived category of a Grothendieck category}, J.\ Eur.\ Math.\ Soc. {\bf 20} (2018), 2607--2641.

\bibitem{Co} L.\ Cohn, \emph{Differential graded categories are $k$-linear stable infinity categories}, arXiv:1308.2587.

\bibitem{Do} M.\ Doni, \emph{$k$-linear Morita theory}, in preparation.

\bibitem{Dr} V.\ Drinfeld, \emph{DG quotients of DG categories}, J.\ Algebra {\bf 272} (2004), 643--691.

\bibitem{GR} F.\ Genovese, J.\ Ramos Gonz\'alez, \emph{A derived Gabriel--Popescu theorem for t-structures via derived injectives}, Int.\ Math.\ Research Notices Volume 2023, Issue 6, 4695--4760.

\bibitem{K2} B.\ Keller, \emph{On differential graded categories}, International Congress of Mathematicians Vol.\ II, Eur.\ Math.\ Soc., Z\"urich (2006), 151--190.

\bibitem{K1} B.\ Keller, \emph{A remark on the generalized smashing conjecture}, Manuscripta Math.\ {\bf 84} (1994), 193--198.

\bibitem{K3} B.\ Keller, \emph{Deriving DG categories}, Ann.\ Sci.\ \'Ecole Norm.\ Sup.\ {\bf 27} (1994), 63--102.

\bibitem{Ko} M.\ Kontsevich, \emph{Homological algebra of Mirror Symmetry}, in: Proceedings of the International Congress of Mathematicians (Zurich, 1994, ed. S.D. Chatterji), Birkhauser, Basel (1995), 120--139.

\bibitem{Kr1} H.\ Krause, \emph{A Brown representability theorem via coherent functors}, Topology {\bf 41} (2002), 853--861.

\bibitem{Lo} A.\ Lorenzin, \emph{Formality and strongly unique enhancements}, arXiv:2204.09527.

\bibitem{LO} V.\ Lunts, D.\ Orlov, \emph{Uniqueness of enhancements for triangulated categories}, J.\ Amer.\ Math.\ Soc.\ {\bf 23} (2010), 853--908.

\bibitem{Nee2} A.\ Neeman, \emph{Bounded t-structures on the category of perfect complexes}, to appear in: Acta Math., https://arxiv.org/abs/2202.08861.

\bibitem{Nee1} A.\ Neeman, \emph{The connection between the K-theory localization theorem of Thomason, Trobaugh and Yao and the smashing subcategories of Bousfield and Ravenel}, Ann.\ Sci.\ \'{E}cole Norm.\ Sup.\ {\bf 25} (1992), 547--566.

\bibitem{Nee3} A.\ Neeman, \emph{Triangulated categories with a single compact generator and a Brown representability theorem}, https://arxiv.org/abs/1804.02240.

\bibitem{Ol} N.\ Olander, \emph{Orlov’s Theorem in the smooth proper case}, J.\ Algebra {\bf 643} (2024), 284--293.

\bibitem{Rickard91} J.\ Rickard, \emph{Derived equivalences as derived functors}, J.\ London Math.\ Soc.\ {\bf 43} (1991), 37--48.

\bibitem{Rickard89b} J.\ Rickard, \emph{Morita theory for derived categories}, J.\ London Math.\ Soc.\ {\bf 39} (1989), 436--456.

\bibitem{R} R.\ Rouquier, \emph{Dimensions of triangulated categories}, J.\ K-theory {\bf 1} (2008), 193--258.

\bibitem{Saorin-Stovicek20} M.\ Saor{\'{\i}}n, J.\ {\v{S}}{\v{t}}ov{\'{\i}}{\v{c}}ek, \emph{$t$-structures with Grothendieck hearts via functor categories}, Selecta Math.\ {\bf 29} (2023), Paper No. 77, 73 pp.

\bibitem{Stack} The Stacks Project Authors, \emph{The Stacks Project}, {\tt http://stacks.math.columbia.edu/}.

\bibitem{To} B.\ To\"en, \emph{The homotopy theory of dg-categories and derived Morita theory}, Invent.\ Math.\ {\bf 167} (2007), 615--667.

\end{thebibliography}
\end{document}